\documentclass[preprint,8pt]{elsarticle}
\usepackage[paperwidth=192mm,paperheight=260mm,top=25mm,bottom=25mm, left=12.5mm, right=12.5mm]{geometry}

\usepackage{graphicx}
\usepackage[font=small,labelfont=bf]{caption}
\usepackage{subcaption}
\usepackage{color}
\usepackage{mdframed}
\usepackage{framed} 
\usepackage{cancel}
\usepackage[normalem]{ulem}

\usepackage{multicol}

\usepackage{tikz}
\usetikzlibrary{tikzmark}
\NewDocumentCommand{\Frame}{}{\Block[tikz={draw,line width= 1.5pt, dashed}]{1-1}{}}

\usepackage{bookmark}
\usepackage{cancel}
\usepackage{multirow}
\usepackage{hhline}
\usepackage{diagbox}
\usepackage{makecell}
\usepackage{slashbox,booktabs}

\usepackage{url}

\usepackage{enumerate}
\usepackage{here}

\usepackage[normalem]{ulem}
\usepackage{ragged2e}
\usepackage{csquotes}
\usepackage{tasks}
\usepackage{paralist}
\usepackage[inline]{enumitem}
\usepackage{tikz-cd}
\usepackage{stmaryrd}

\usepackage{afterpage}
\usepackage{lineno,hyperref}
\modulolinenumbers[5]

\usepackage{amsmath,amssymb,amsthm,mathtools,mathrsfs,epsfig,amsfonts,wasysym}
\usepackage{MnSymbol}

\usepackage[normalem]{ulem}
\usepackage{ragged2e}
\usepackage{csquotes}
\usepackage{tasks}
\usepackage{paralist}
\usepackage[inline]{enumitem}
\usepackage{tikz-cd}
\usepackage{stmaryrd}

\usepackage[titletoc,title]{appendix}

\newcommand{\e}{\mathrm{e}}
\newcommand{\de}{\mathrm{d}}

\usepackage{url,hyperref}
\hypersetup{colorlinks=true,linkcolor=blue,filecolor=magenta,urlcolor=cyan}
\usepackage{nicematrix}
\usepackage{tikz}
\usetikzlibrary{fit}

\usepackage[all]{xy}

\usepackage{booktabs}
\usepackage{tabularray}
\UseTblrLibrary{diagbox}
\usepackage{colortbl}

\newtheorem{definition}{Defintion}

\newtheorem{proposition}{Proposition}

\usepackage{natbib}
\setcitestyle{authoryear,open={(},close={)}}

\usepackage{pifont}
\newcommand{\xmark}{\ding{55}}
\usepackage{float}

\usepackage{orcidlink}

\definecolor{rose}{rgb}{1.0, 0.0, 0.50}

\usepackage{amssymb}

\journal{Infectious Disease Modelling}

\begin{document}

\begin{frontmatter}



\title{A novel comparison framework for epidemiological strategies applied to age-based restrictions versus horizontal lockdowns}


\author[mymainaddress1]{Vasiliki Bitsouni \corref{mycorrespondingauthor}\negthickspace \orcidlink{0000-0002-0684-0583}}
\cortext[mycorrespondingauthor]{Corresponding author}
\ead{vbitsouni@math.upatras.gr}
\ead[url]{https://www.math.upatras.gr/el/people/vbitsouni}
\author[mymainaddress2,mymainaddress3]{Nikolaos Gialelis \!\orcidlink{0000-0002-6465-7242}}
\ead{ngialelis@math.uoa.gr}
\ead[url]{http://users.uoa.gr/~ngialelis}
\author[mymainaddress1]{Vasilis Tsilidis  \!\orcidlink{0000-0001-5868-4984}}
\ead{vtsilidis@upatras.gr}
\ead[url]{https://tsilidisv.github.io/}

\address[mymainaddress1]{Department of Mathematics, University of Patras, GR-26504 Rio Patras, Greece}

\address[mymainaddress2]{Department of Mathematics, National and Kapodistrian University of Athens, GR-15784 Athens, Greece}

\address[mymainaddress3]{School of Medicine, National and Kapodistrian University of Athens, GR-11527 Athens, Greece}


\begin{abstract}
During an epidemic, such as the COVID-19 pandemic, policy-makers are faced with the decision of implementing effective, yet socioeconomically costly intervention strategies, such as school and workplace closure, physical distancing, etc. In this study, we propose a rigorous definition of epidemiological strategies. In addition, we develop a scheme for comparing certain epidemiological strategies, with the goal of providing policy-makers with a tool for their systematic comparison. Then, we put the suggested scheme to the test by employing an age-based epidemiological compartment model introduced in \cite{bgt2023SVEAIRmodel},  coupled with data from the literature, in order to compare the effectiveness of age-based and horizontal interventions. In general, our findings suggest that these two are comparable, mainly at a low or medium level of intensity.

\end{abstract}




\end{frontmatter}



\section{Introduction}
\label{intro}

The recent COVID-19 pandemic brought to the fore the disastrous for the economy consequences of horizontal lockdowns. Economically costly horizontal measures during the COVID-19 pandemic have been the closure of workplaces and schools, the cancellation of public events and general stay-at-home restrictions (see \cite{brodeur2021literature,chen2021epidemiological,Deb2021,ourworldindata} and many references therein). 

This fact highlights the need for a more sophisticated managing of epidemiological crises. In this context, many countries, especially after the spasmodic first response, have looked for more flexible intervention policies. Multiple combinations of interventions were deployed by policy-makers in order to combat the spread of SARS-CoV-2 and minimize their impact on the economy \citep{asahi2021effect,karatayev2020local,perra2021non}. 

Finding ways to intervene in the natural progression of disease spreading, has been a hot topic in the scientific community. Models have been proposed, for a wide range of diseases, investigating various non-pharmaceutical interventions \citep{demers2023relationship, adegbite2023mathematical, verma2020capacity, zakary2017new, bhadauria2023studying, vatcheva2021social, brethouwer2021stay, amaku2021modelling, saha2022effect}, vaccination \citep{paton2023evaluation, anupong2023modeling, owusu2023equitable, gan2024need, thongtha2022optimal, abell2023understanding} and treatment \citep{ZAMAN200943, beraud2022remdesivir} strategies, as well as various combinations of the aforementioned interventions \citep{apenteng2020impact, lamba2024impact}. Despite the success of the foregoing studies, the lack of a mathematically rigorous definition of epidemiological strategies is apparent.

Moreover, age-based interventions have been discussed as a theoretical alternative to horizontal lockdowns. However, they have also raised ethical concerns with regard to ageism \citep{van2020age, spaccatini2022you, motorniak2023reelin}.

To our knowledge, the investigation of age-based interventions has been limited in terms of modeling. The authors of \cite{acemoglu2021optimal} proposed a multigroup SIR model, with the intent of studying age-based lockdowns. In \cite{kirwin2021net}, the authors study the prioritization of vaccination to selected target groups. 

In the present study, we:
\begin{itemize}
    \item[--] give a rigorous definition of the notion of epidemiological strategies
    \item[--] propose a framework for systematically comparing certain epidemiological strategies
    \item[--] utilize the aforementioned scheme to compare the effectiveness of age-based interventions when compared to horizontal lockdowns, in the case of the SARS-CoV-2 pandemic.
\end{itemize}

This study is organized as follows. In \hyperref[sec:versus]{\S \ref*{sec:versus}}, we introduce the notion of an (epidemiological) strategy, along with its potential gradations, and we present a framework for comparing the effectiveness of certain strategies.  In \hyperref[sec:VS]{\S \ref*{sec:VS}}, we contrast the impact of a horizontal lockdown with varying levels of intensity, with certain age-based countermeasures that have a similar epidemiological effect, but less of an influence on society and, consequently, the economy. In \hyperref[sec:CD]{\S \ref*{sec:CD}}, we conclude with a summary and discussion of the results.


\section{A framework for comparing the effectiveness of different strategies}
\label{sec:versus}

Let us divide a population into two classes, the infectious, $\mathcal{I}$, and the non-infectious, $\mathcal{I}^c$. Each of these classes can be divided to further sub-compartments, e.g., $A\in\mathcal{I}$ and $B\in\mathcal{I}^c$. 

\begin{figure}[h]
    \centering
    \includegraphics[width=.5\textwidth]{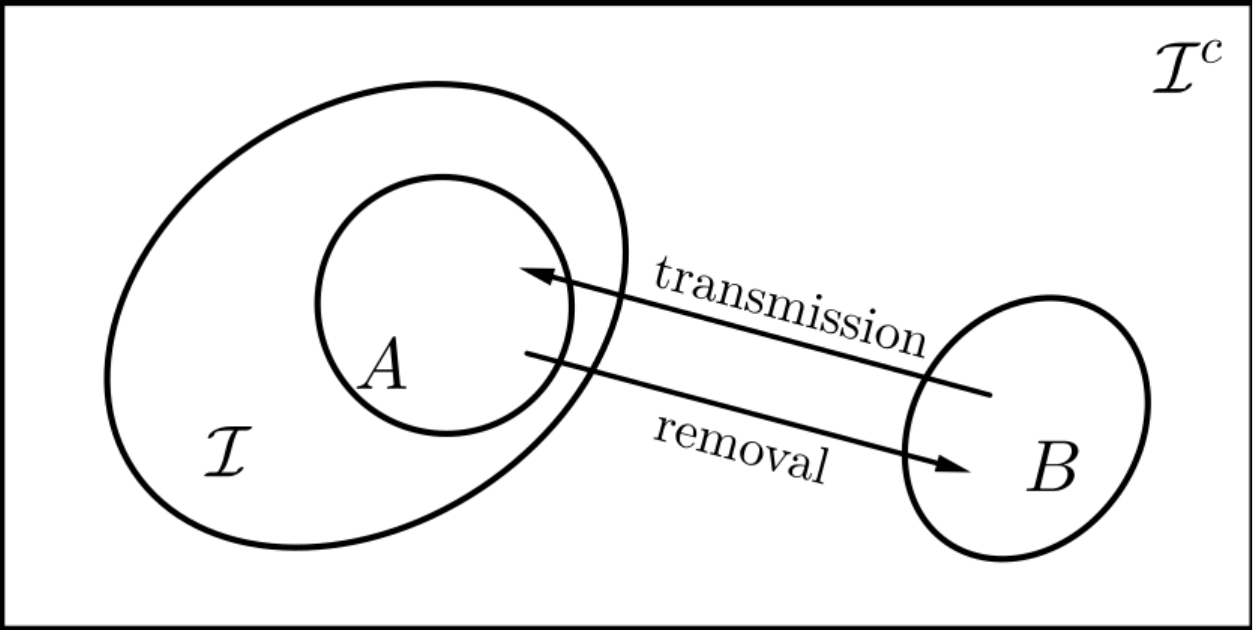}
    \caption{Flows between the classes of infectious, $\mathcal{I}$, and non-infectious, $\mathcal{I}^c$, individuals of a population.}
    \label{fig:1}
\end{figure}

The \textit{transmission rate} from compartment $B$ to compartment $A$ is defined as 
\begin{equation}
    \beta_{B\to A}\coloneqq \frac{c_B \cdot\varpi_{B\to A}}{N} \; ,
    \label{beta-def}
\end{equation}
where $c_B$ is the average number of close contacts of an individual belonging in $B$ with other individuals, $\varpi_{B\to A}$ is the probability of a contact to be effective in turning an individual of compartment $B$ to an individual of compartment $A$, and $$N\coloneqq \mathcal{I}+\mathcal{I}^c$$ is the total number of the population. The \textit{removal rate} from compartment $A$ to compartment $B$ is defined as 
\begin{equation}
    \gamma_{A\to B}\coloneqq \frac{1}{P_{A\to B}} \; ,
    \label{gamma-def}
\end{equation}
where $P_{A\to B}$ is the average period an individual spends on compartment $A$ before moving into compartment $B$. A diagram for the above definitions is shown in \hyperref[fig:1]{Figure \ref*{fig:1}}. 

These parameters, probably among others depending on the model (for instance, the model employed later on in the present study comprises two types of transmission rates and two types of removal rates, among nine other parameters), are involved into the formulation of an epidemiological model that describes an epidemiological problem under study. However, these parameters are special, because interventions by external factors acting for the control of the studied epidemiological phenomenon (e.g., policy makers), can be described as changes in their values.

Now, enumerating all the different transmission and removal rates of a particular model, i.e., $\beta_1,\dots,\beta_{n_1}$ and $\gamma_1,\dots,\gamma_{n_2}$, respectively, we can write $$\pmb{\beta}={\left(\beta_i\right)}_{i=1}^{n_1}\text{ and }\pmb{\gamma}={\left(\gamma_{i}\right)}_{i=1}^{n_2}.$$ Throughout the present section we assume a well-posed global (with respect to time) epidemiological compartmental problem, $$\mathscr{P}=\mathscr{P}{\left(\mathscr{M}\right)},$$ which is described by a (differential/difference equations, agent-based, etc.) model, $$\mathscr{M}=\mathscr{M}{\left(\pmb{x},\left(\pmb{\beta}{\left(\pmb{x}\right)},\pmb{\gamma}{\left(\pmb{x}\right)}\right),\pmb{\delta}{\left(\pmb{x}\right)}\right)},$$ where $$\pmb{x}={\left(x_i\right)}_{i=1}^m \in\mathcal{X}\subseteq\mathbb{R}^m$$ is the vector of the independent variables, $$\left(\pmb{\beta},\pmb{\gamma}\right)\in F{\left(\mathcal{X};\mathcal{P}_{\mathrm{tr},\mathrm{r}}\subseteq\,\mathbb{R}^{n_1}\times\mathbb{R}^{n_2}\right)}=\left\{\pmb{f}\colon\,\mathcal{X}\to\mathcal{P}_{\mathrm{tr},\mathrm{r}}\right\}$$ is the vector-valued function of the parameters of interest of $\mathscr{M}$ and $$\pmb{\delta}\in F{\left(\mathcal{X};\mathcal{P}_{\mathrm{other}}\subseteq\mathbb{R}^{n_3}\right)}$$ is the vector-valued function of the rest of parameters of $\mathscr{M}$. 

\subsection{Strategies and substrategies}
\label{strategies}

We begin by introducing the concept of a \textit{strategy} of $\mathscr{P}$, which is of pivotal importance for the following analysis. In the epidemiology framework, a strategy can be considered as the mathematical description of a set of epidemiological interventions made by potential external factors, such as policy makers, experts etc., in order to restrict the epidemiological phenomenon. These interventions consist of first fixing a reference value, $\left(\pmb{\beta}_0,\pmb{\gamma}_0\right)\in F{\left(\mathcal{X};\mathcal{P}_{\mathrm{tr},\mathrm{r}}\right)}$, for the parameters chosen, and then scaling each element of the set in terms of the fixed value.

\begin{definition}[strategy \& strategic scale of an element]
\label{strdef1}
Let $\left(\pmb{\beta}_0,\pmb{\gamma}_0\right)\in F{\left(\mathcal{X};\mathcal{P}_{\mathrm{tr},\mathrm{r}}\right)}$. 
\begin{enumerate}
\item A set $S=S{\left(\pmb{\beta}_0,\pmb{\gamma}_0\right)}\subseteq F{\left(\mathcal{X};\mathcal{P}_{\mathrm{tr},\mathrm{r}}\right)}$ is called strategy (of $\mathscr{P}$) with respect to $\left(\pmb{\beta}_0,\pmb{\gamma}_0\right)$ iff $$\forall \pmb{y}\in S\text{ }\exists\pmb{r}=\pmb{r}{\left(\,\cdot\,;\left(\pmb{\beta}_0,\pmb{\gamma}_0\right),\pmb{y}\right)}\in F{\left(\mathcal{X};\mathbb{R}^{n_1+n_2}\right)}\text{ s.t. }\pmb{y}=\pmb{r}\odot\left(\pmb{\beta}_0,\pmb{\gamma}_0\right),$$ where $\odot$ stands for the Hadamard product. 
\item Let 
\begin{enumerate}[label=\roman*.]
\item $S=S{\left(\pmb{\beta}_0,\pmb{\gamma}_0\right)}$ be a respective  strategy and 
\item $\pmb{y}\in S$.
\end{enumerate}
A function $\pmb{r}\in F{\left(\mathcal{X};\mathbb{R}^{n_1+n_2}\right)}$ as in $1$. is called strategic scale of $\pmb{y}$. 
\end{enumerate}
\end{definition} 

We observe that every subset of a strategy is a strategy itself, as it is referred in the following result, the elementary proof of which is omitted. 
\begin{proposition}
\label{strpr1} 
Let 
\begin{enumerate}
\item $\left(\pmb{\beta}_0,\pmb{\gamma}_0\right)\in F{\left(\mathcal{X};\mathcal{P}_{\mathrm{tr},\mathrm{r}}\right)}$, 
\item $S=S{\left(\pmb{\beta}_0,\pmb{\gamma}_0\right)}$ be a respective strategy and 
\item $S_0\subseteq S$.
\end{enumerate}
Then $S_0$ is a strategy with respect to $\left(\pmb{\beta}_0,\pmb{\gamma}_0\right)$. 
\end{proposition}

In view of \hyperref[strpr1]{Proposition \ref*{strpr1}}, we give the definition of a substrategy of a given strategy. In the epidemiology framework, a \textit{substrategy} can be considered as the mathematical description of a subset of a given set of epidemiological interventions. 
\begin{definition}[substrategy]
\label{strdef2} 
Let 
\begin{enumerate}[label=\roman*.]
\item $\left(\pmb{\beta}_0,\pmb{\gamma}_0\right)\in F{\left(\mathcal{X};\mathcal{P}_{\mathrm{tr},\mathrm{r}}\right)}$, 
\item $S=S{\left(\pmb{\beta}_0,\pmb{\gamma}_0\right)}$ be a respective strategy and 
\item $S_0\subseteq S$. 
\end{enumerate} 
We call $S_0$ a substrategy of $S$. 
\end{definition}

In fact, we can define a substrategy by setting limitations to the choice of a strategic scale of each of its elements. Below, we name certain such examples. 
\begin{definition}[horizontal and $x_i$-based strategy]
\label{strdef3} 
Let 
\begin{enumerate}[label=\roman*.]
\item $\left(\pmb{\beta}_0,\pmb{\gamma}_0\right)\in F{\left(\mathcal{X};\mathcal{P}_{\mathrm{tr},\mathrm{r}}\right)}$, 
\item $S=S{\left(\pmb{\beta}_0,\pmb{\gamma}_0\right)}$ be a respective strategy and 
\item $S_0\subseteq S$. 
\end{enumerate} 
We name the following substrategies. 
\begin{enumerate} 
\item Let $i\in\left\{1,\dots,m\right\}$. $S_0$ is called horizontal with respect to $x_i$ iff $$\pmb{r}{\left(\pmb{x};\left(\pmb{\beta}_0,\pmb{\gamma}_0\right),\pmb{y}\right)}=\pmb{r}{\left(x_1,\dots,x_{i-1},x_{i+1},\dots,x_m;\left(\pmb{\beta}_0,\pmb{\gamma}_0\right),\pmb{y}\right)},\text{ }\forall\pmb{x}\in\mathcal{X},\text{ }\forall\pmb{y}\in S_0,$$ i.e., $\forall\pmb{y}\in S_0$ a respective strategic scale is independent of $x_i$, otherwise we call it $x_i$-based. 
\item $S_0$ is called horizontal, iff it is horizontal with respect to $x_i$, $\forall i\in\left\{1,\dots,m\right\}$. 
\end{enumerate}
\end{definition}
In the epidemiology framework, a $x_i$-based substrategy can be considered as the mathematical description of a subset of epidemiological interventions, which targets a certain group of a population partitioned with respect to $x_i$ variable. 

We also observe that every union of strategies is a strategy itself, as it is referred in the following elementary result. 
\begin{proposition}
\label{strpr2} 
Let 
\begin{enumerate}
\item $\left(\pmb{\beta}_0,\pmb{\gamma}_0\right)\in F{\left(\mathcal{X};\mathcal{P}_{\mathrm{tr},\mathrm{r}}\right)}$ and 
\item ${\left\{S_j=S_j{\left(\pmb{\beta}_0,\pmb{\gamma}_0\right)}\right\}}_{j\in\mathcal{J}}$ be a family of respective strategies. 
\end{enumerate} 
Then $\bigcup\limits_{j\in\mathcal{J}}{S_j}$ is a strategy with respect to $\left(\pmb{\beta}_0,\pmb{\gamma}_0\right)$. 
\end{proposition}

In view of \hyperref[strpr2]{Proposition \ref*{strpr2}}, we give the following definition. 
\begin{definition}[the largest strategy]
\label{strdef4} 
Let 
\begin{enumerate}
\item $\left(\pmb{\beta}_0,\pmb{\gamma}_0\right)\in F{\left(\mathcal{X};\mathcal{P}_{\mathrm{tr},\mathrm{r}}\right)}$ and 
\item $\mathcal{S}$ be the family of all the respective strategies. 
\end{enumerate}
We call $$\hat{S}=\hat{S}{\left(\pmb{\beta}_0,\pmb{\gamma}_0\right)}=\bigcup\limits_{S\in\mathcal{S}}{S}$$ the largest strategy with respect to $\left(\pmb{\beta}_0,\pmb{\gamma}_0\right)$. 
\end{definition}

Of course, whatever result holds for the largest strategy also holds for an abstract strategy, like the following direct one.  
\begin{proposition}
\label{strpr3} 
Let $\left(\pmb{\beta}_0,\pmb{\gamma}_0\right)\in F{\left(\mathcal{X};\mathcal{P}_{\mathrm{tr},\mathrm{r}}\right)}$. If 
\begin{equation}
\label{uniqcond}
\left(\underbrace{0,\dots,0}_{\#\text{ }n_1+n_2}\right)\nin\left(\pmb{\beta}_0{\left(\mathcal{X}\right)},\pmb{\gamma}_0{\left(\mathcal{X}\right)}\right),
\end{equation}
then $$\forall\pmb{y}\in \hat{S}\text{ }\exists!\text{ strategic scale of }\pmb{y}.$$ 
\end{proposition}

Under the light of \hyperref[strpr3]{Proposition \ref*{strpr3}}, the next notion is well-defined. 
\begin{definition}[strategic scale of a strategy]
\label{strdef5}
Let 
\begin{enumerate}
\item $\left(\pmb{\beta}_0,\pmb{\gamma}_0\right)\in F{\left(\mathcal{X};\mathcal{P}_{\mathrm{tr},\mathrm{r}}\right)}$ satisfy \eqref{uniqcond} and 
\item $S=S{\left(\pmb{\beta}_0,\pmb{\gamma}_0\right)}$ be a respective strategy. 
\end{enumerate}
We call the function $$S\ni\pmb{y}\mapsto\pmb{r}{\left(\,\cdot\,;\left(\pmb{\beta}_0,\pmb{\gamma}_0\right),\pmb{y}\right)}\in F{\left(\mathcal{X};\mathbb{R}^{n_1+n_2}\right)}$$ the strategic scale of $S$. 
\end{definition} 

We can then easily deduce the following result. 
\begin{proposition}
\label{strpr4} 
Let 
\begin{enumerate}[label=\roman*.]
\item $\left(\pmb{\beta}_0,\pmb{\gamma}_0\right)\in F{\left(\mathcal{X};\mathcal{P}_{\mathrm{tr},\mathrm{r}}\right)}$ satisfy \eqref{uniqcond}, 
\item $S=S{\left(\pmb{\beta}_0,\pmb{\gamma}_0\right)}$ be a respective strategy,  
\item $\mathfrak{P}_j$, $\forall j\in\left\{1,2\right\}$ be mathematical statements with respect to the strategic scale of $S$ such that $\mathfrak{P}_1\Rightarrow\mathfrak{P}_2$ and 
\item $S_j\subseteq S$, $\forall j\in\left\{1,2,3\right\}$, be substrategies of $S$ such that $$S_j=\left\{\pmb{y}\in S\,\big|\,\mathfrak{P}_i{\left(\pmb{r}{\left(\,\cdot\,;\left(\pmb{\beta}_0,\pmb{\gamma}_0\right),\pmb{y}\right)}\right)}\right\},\text{ }\forall j\in\left\{1,2\right\}\text{ \& }S_3=\left\{\pmb{y}\in S\,\big|\,\neg\mathfrak{P}_1{\left(\pmb{r}{\left(\,\cdot\,;\left(\pmb{\beta}_0,\pmb{\gamma}_0\right),\pmb{y}\right)}\right)}\right\}.$$
\end{enumerate}
Then 
\begin{enumerate}
\item $S_1\subseteq S_2$ and 
\item $S_3=S\setminus S_1$. 
\end{enumerate}
\end{proposition} 

For example, for a given $\left(\pmb{\beta}_0,\pmb{\gamma}_0\right)\in F{\left(\mathcal{X};\mathcal{P}_{\mathrm{tr},\mathrm{r}}\right)}$ that satisfies \eqref{uniqcond} and a given $i\in\left\{1,\dots,m\right\}$, the horizontal with respect to $x_i$ substrategy of $\hat{S}=\hat{S}{\left(\pmb{\beta}_0,\pmb{\gamma}_0\right)}$ is the set-theoretic complement with respect to $\hat{S}$ of the $x_i$-based substrategy of $\hat{S}$. The scope of the present paper can be now stated as the comparison of the above substrategies for $x_i$ being the age of an individual of a population. 

\subsection{Comparison of strategies}
\label{comparison} 

Here we introduce a scheme for the comparison of strategies, for which we need some preliminary notions, such as the basic reproductive number and the graded strategies. 

\subsubsection{\texorpdfstring{$\mathcal{R}_0$}{R0}: the measure of comparison}
\label{r0coin} 

An important epidemiological notion studied and used extensively in the epidemiological literature is the basic reproductive number, $\mathcal{R}_0$, which is defined as the average number of infectious cases directly generated by one such case in a population where all individuals are susceptible to an infection. For every mathematical model, that describes a problem under study, corresponds a respective $\mathcal{R}_0$, which can be calculated with several ways, such as with the next-generation method or the existence of the endemic equilibrium \citep{diekmann2000mathematical}. 

In general, $\mathcal{R}_0$ depends on both the independent variables and the parameters of a model, therefore it is considered as a function defined as 
\begin{align*}
\mathcal{R}_0\colon\,\mathcal{X}\times F{\left(\mathcal{X};\mathcal{P}_{\mathrm{tr}},\mathrm{r}\right)}\times F{
\left(\mathcal{X};\mathcal{P}_{\mathrm{other}}\right)}&\to\left(0,\infty\right)\\
\left(\pmb{x},\left(\pmb{\beta},\pmb{\gamma}\right),\pmb{\delta}\right)&\mapsto\mathcal{R}_0{\left(\pmb{x},\left(\pmb{\beta},\pmb{\gamma}\right),\pmb{\delta}\right)}.
\end{align*}
Only for the sake of brevity and compactness of the exposition, in the present paper we assume that it is independent of $\pmb{x}$, that is 
\begin{align*}
\mathcal{R}_0\colon\,F{\left(\mathcal{X};\mathcal{P}_{\mathrm{tr}},\mathrm{r}\right)}\times F{
\left(\mathcal{X};\mathcal{P}_{\mathrm{other}}\right)}&\to\left(0,\infty\right)\\
\left(\left(\pmb{\beta},\pmb{\gamma}\right),\pmb{\delta}\right)&\mapsto\mathcal{R}_0{\left(\left(\pmb{\beta},\pmb{\gamma}\right),\pmb{\delta}\right)}.
\end{align*}

In the proposed scheme, we check how one strategy measures against another of a special kind, via the calculation of the respective values of $\mathcal{R}_0$. That special kind of strategies is described below. 

\subsubsection{Gradable and graded strategies}
\label{grgrstr} 

The notion of the graded strategies is the crux of the proposed scheme. But before its introduction, we first need the following one. 
\begin{definition}[gradable strategy]
\label{grdef1} 
Let 
\begin{enumerate}
\item $\left(\pmb{\beta}_0,\pmb{\gamma}_0\right)\in F{\left(\mathcal{X};\mathcal{P}_{\mathrm{tr},\mathrm{r}}\right)}$ satisfy \eqref{uniqcond} and 
\item $S=S{\left(\pmb{\beta}_0,\pmb{\gamma}_0\right)}$ be a respective strategy. 
\end{enumerate}
We call $S$ gradable iff $\forall\pmb{\delta}\in F{\left(\mathcal{X};\mathcal{P}_{\mathrm{other}}\right)}$ the function $\left.\mathcal{R}_0\right|_S{\left(\,\cdot\,,\pmb{\delta}\right)}$ is injective.
\end{definition}

Since $\mathcal{R}_0{\left(S,\pmb{\delta}\right)}\subseteq\left(0,\infty\right)$, $\forall\pmb{\delta}\in F{\left(\mathcal{X};\mathcal{P}_{\mathrm{other}}\right)}$, we can arrange any family of pairwise distinct elements of such a set in a strictly ascending order when $S$ is gradable, hence the following notion is well-defined. 
\begin{definition}[graded strategy]
\label{grdef2} 
Let 
\begin{enumerate}
\item $\left(\pmb{\beta}_0,\pmb{\gamma}_0\right)\in F{\left(\mathcal{X};\mathcal{P}_{\mathrm{tr},\mathrm{r}}\right)}$ satisfy \eqref{uniqcond}, 
\item $S=S{\left(\pmb{\beta}_0,\pmb{\gamma}_0\right)}$ be a respective gradable strategy, 
\item $\pmb{\delta}\in F{\left(\mathcal{X};\mathcal{P}_{\mathrm{other}}\right)}$ and 
\item ${\left\{\pmb{y}_i\right\}}_{i=1}^k\subseteq S$ be a family of pairwise distinct elements of $S$, such that $$\underbrace{\mathcal{R}_0{\left(\pmb{y}_1,\pmb{\delta}\right)}}_{\eqqcolon G_1}<\dots<\underbrace{\mathcal{R}_0{\left(\pmb{y}_k,\pmb{\delta}\right)}}_{\eqqcolon G_k}.$$  
\end{enumerate}
We call the pair $\left(S,\pmb{G}={\left(G_i\right)}_{i=1}^k\right)$ a graded strategy, while $\pmb{G}$ is called a gradation of $S$ and each of the $G_1,\dots,G_k$ is called a grade of $\pmb{G}$. 
\end{definition}

We note that the gradation of a graded strategy is a matter of choice. In what follows, for the sake of simplicity, we write $S$ instead of $\left(S,\pmb{G}\right)$ for a graded strategy. 

\subsubsection{Comparison table and coverage}
\label{epcov} 

With the above toolbox at hand, we propose a scheme for the comparison of two strategies, only when one of them is graded. In addition, the scheme allows us to include many substrategies of the other strategy. Below, we present the steps required for the utilization of the proposed scheme, which is governed by the construction of the respective \textit{comparison table} and analyzed in terms of \textit{epidemiological} and \textit{social coverage}.

\paragraph{Construction of the comparison table} 

\definecolor{DustyGray}{rgb}{0.58,0.58,0.58}

\begin{enumerate}
\item Placing of the grades $G_1,\dots,G_k$ of the gradation $\pmb{G}={\left(G_i\right)}_{i=1}^k$ of a given graded strategy $S_1$, in increasing order, to the top row:
\begin{center}
    \begin{tabular}[t]{cllllll}
    \toprule
    \diagbox{\cellcolor{DustyGray}}{$S_1$} & $G_1$ & $G_2$ & $\cdots$ & $G_k$  \\ \cmidrule(r){1-1} \cmidrule(rl){2-5} 
     \rowcolor{DustyGray}
                  &       &      &         &        \\
    \rowcolor{DustyGray}              &       &       &          &       \\
    \rowcolor{DustyGray}               &       &       &          &        \\
    \rowcolor{DustyGray}             &       &       &          &        \\
    \bottomrule
    \end{tabular}
\end{center}
\item Placing the under study substrategies ${\left\{S_{2_i}\right\}}_{i=1}^\ell$ of a second strategy $S_2$ to the left of the table, with the intent of comparing them against the first graded strategy.
\begin{center}
    \begin{tabular}[t]{cllllll}
    \toprule
    \diagbox{$S_2$}{$S_1$} & $G_1$ & $G_2$ & $\cdots$ & $G_k$  \\ \cmidrule(r){1-1} \cmidrule(rl){2-5} 
    $S_{2_1}$       &   \cellcolor{DustyGray}     &    \cellcolor{DustyGray}  &     \cellcolor{DustyGray}    &    \cellcolor{DustyGray}   \\
    $S_{2_2}$              &    \cellcolor{DustyGray}   &   \cellcolor{DustyGray}    &     \cellcolor{DustyGray}     & \cellcolor{DustyGray}   \\
    $\vdots$              &   \cellcolor{DustyGray}    &  \cellcolor{DustyGray}     &     \cellcolor{DustyGray}     &   \cellcolor{DustyGray}  \\
     $S_{2_\ell}$           &   \cellcolor{DustyGray}    &   \cellcolor{DustyGray}    &   \cellcolor{DustyGray}       & \cellcolor{DustyGray} \\
    \bottomrule
    \end{tabular}
\end{center}
\item Populating the comparison table with $\star$, where
\begin{equation*}
    \star_{ij} = \begin{cases}
                     \checkmark & \text{ if the $\mathcal{R}_0$ of $S_{2i}$ is greater than or equal to $G_j$ } \\
                     \text{\xmark}     & \text{ otherwise}\;, 
                 \end{cases}\text{ }\forall\left(i,j\right)\in\left\{1,\dots,\ell\right\}\times\left\{1,\dots,k\right\}\;.
\end{equation*}
\begin{center}
    \begin{tabular}[t]{ccccc}
    \toprule
    \diagbox{$S_2$}{$S_1$} & $G_1$ & $G_2$ & $\cdots$ & $G_k$   \\ \cmidrule(r){1-1} \cmidrule(rl){2-5} 
    $S_{2_1}$       & $\star_{11}$     &   $\star_{12}$  &     $\cdots$    &   $\star_{1k}$    \\
    $S_{2_2}$              &     $\star_{21}$     &   $\star_{22}$  &     $\cdots$    &   $\star_{2k}$   \\
    $\vdots$              &     $\vdots$     &   $\vdots$  &     $\ddots $    &   $\vdots$     \\
     $S_{2_\ell}$           &     $\star_{\ell1}$     &   $\star_{\ell2}$  &     $\cdots$    &   $\star_{\ell k}$     \\
    \bottomrule
    \end{tabular}
\end{center}
\end{enumerate}

Next, we present two ways to read the comparison table for extracting useful information. 

\paragraph{Social overview of the comparison table: epidemiological coverage} Here we compare $S_2$ to $S_1$. In particular, for every fixed substrategy of $S_2$ (social overview), we check how good of an alternative it is, compared to $S_1$ (epidemiological coverage).

\begin{enumerate}[label=\arabic*$_\rightarrow$.]
\setcounter{enumi}{3}
\item Calculating the epidemiological coverage of the gradation $\pmb{G}$ of $S_1$ by each substrategy of $S_2$, by calculating the average number of $\checkmark$ in each row.
\begin{center}
\begin{NiceTabular}{cccccc}[]
\CodeBefore
\cellcolor[HTML]{949494}{6-6}
\Body
\toprule 
\diagbox{$S_2$}{$S_1$} & $G_1$   & $G_2$   & $\cdots$ & $G_k$   & \makecell{Epidemiological \\coverage ($\cdot100\%$)}  \\
\cmidrule(r){1-1} \cmidrule(rl){2-5} \cmidrule(rl){6-6}            
$S_{2_1}$              & $\star_{11}$     &   $\star_{12}$  &     $\cdots$    &   $\star_{1k}$ & $\frac{\#{\left\{\star_{1 j}=\checkmark\right\}}_{j=1}^k}{k}$           \\
 $S_{2_2}$              & $\star_{21}$     &   $\star_{22}$  &     $\cdots$    &   $\star_{2k}$ &  $\frac{\#{\left\{\star_{2 j}=\checkmark\right\}}_{j=1}^k}{k}$        \\
$\vdots$               & $\vdots$ & $\vdots$ & $\ddots $  & $\vdots$ & $\vdots$                \\
 $S_{2_\ell}$              & $\star_{\ell1}$     &   $\star_{\ell2}$  &     $\cdots$    &   $\star_{\ell k} $ & $\frac{\#{\left\{\star_{\ell j}=\checkmark\right\}}_{j=1}^k}{k}$               \\   \specialrule{.8pt}{2pt}{3pt}
& & & & & \Frame$\cellcolor{DustyGray}\phantom{\frac{\#{\left\{\star_{i j}=\checkmark\right\}}_{i,j=1}^{\ell,k}}{\ell \cdot k}}$  \\

\end{NiceTabular}
\end{center}
\item Calculating the total coverage of the gradation $\pmb{G}$ of $S_1$ by the whole $S_2$, through the average value of all epidemiological coverages. 
\begin{center}
\begin{NiceTabular}[t]{cccccc}
\toprule 
\diagbox{$S_2$}{$S_1$} & $G_1$   & $G_2$   & $\cdots$ & $G_k$   & \makecell{Epidemiological \\coverage ($\cdot100\%$)}  \\
\cmidrule(r){1-1} \cmidrule(rl){2-5} \cmidrule(rl){6-6}            
$S_{2_1}$              & $\star_{11}$     &   $\star_{12}$  &     $\cdots$    &   $\star_{1k}$ & $\frac{\#{\left\{\star_{1 j}=\checkmark\right\}}_{j=1}^k}{k}$           \\
$S_{2_2}$              & $\star_{21}$     &   $\star_{22}$  &     $\cdots$    &   $\star_{2k}$ &  $\frac{\#{\left\{\star_{2 j}=\checkmark\right\}}_{j=1}^k}{k}$        \\
$\vdots$               & $\vdots$ & $\vdots$ & $\ddots $  & $\vdots$ & $\vdots$                \\
$S_{2_\ell}$              & $\star_{\ell1}$     &   $\star_{\ell2}$  &     $\cdots$    &   $\star_{\ell k} $ & $\frac{\#{\left\{\star_{\ell j}=\checkmark\right\}}_{j=1}^k}{k}$               \\   \specialrule{.8pt}{2pt}{3pt}
& & & & & \Frame$\frac{\#{\left\{\star_{i j}=\checkmark\right\}}_{i,j=1}^{\ell,k}}{\ell \cdot k}$  \\
\end{NiceTabular}
\end{center}
\end{enumerate}

The takeaway of the above analysis is that if the total (epidemiological) coverage of $\pmb{G}$ by the (respective sub-)strategy $S_2$ ($S_{2_i}$, for $i\in\left\{1,\dots,\ell\right\}$) is satisfying, then $S_2$ ($S_{2_i}$) could be considered as an alternative to $S_1$. We note that the quantification of the term \textquote{satisfying} is subjective.

\paragraph{Epidemiological overview of the comparison table: social coverage} Here we compare $S_1$ to $S_2$. In particular, for every fixed grade of $\pmb{G}$ (epidemiological overview), we check how well it can be covered by $S_2$ (social coverage). 

\begin{enumerate}[label=\arabic*$_\downarrow$.]
\setcounter{enumi}{3}
\item Calculating the social coverage of $S_2$ by each grade of $\pmb{G}$ of $S_1$, by calculating the average number of $\checkmark$ in each column.
\begin{center}
\begin{NiceTabular}[t]{cccccc}
\CodeBefore
\cellcolor[HTML]{949494}{6-6}
\Body
\cmidrule[.08em]{1-5}  
\Block{1-1}
  {
    \diagbox
      {\rule[-2mm]{0pt}{2mm}$S_2$}
      {\rule{0pt}{4mm}$S_1$ }
  }  & $G_1$   & $G_2$   & $\cdots$ & $G_k$  \\
\cmidrule(lr){1-1} \cmidrule(rl){2-5}
$S_{2_1}$              & $\star_{11}$     &   $\star_{12}$  &     $\cdots$    &   $\star_{1k}$ &       \\
$S_{2_2}$              & $\star_{21}$     &   $\star_{22}$  &     $\cdots$    &   $\star_{2k}$ &                   \\
$\vdots$               & $\vdots$ & $\vdots$ & $\ddots $  & $\vdots$ &                \\
$S_{2_\ell}$              & $\star_{\ell1}$     &   $\star_{\ell2}$  &     $\cdots$    & $\star_{\ell k}$        \\   
\cmidrule(lr){1-1} \cmidrule(lr){2-5}
Social coverage $(\cdot100\%)$ & $\frac{\#{\left\{\star_{i 1}=\checkmark\right\}}_{i=1}^\ell}{\ell}$ & $\frac{\#{\left\{\star_{i 2}=\checkmark\right\}}_{i=1}^\ell}{\ell}$ & $\cdots$ & $\frac{\#{\left\{\star_{i k}=\checkmark\right\}}_{i=1}^\ell}{\ell}$ & \Frame \cellcolor{DustyGray}\phantom{$\frac{\#{\left\{\star_{i j}=\checkmark\right\}}_{i,j=1}^{\ell,k}}{\ell \cdot k}$} \\
\cmidrule[.08em]{1-5} 
\end{NiceTabular}
\end{center}
\item Calculating the total coverage of 
$S_2$ by the whole $\pmb{G}$ of $S_1$, through the average value of all social coverages. 
\begin{center}
\begin{NiceTabular}[t]{cccccc}
\cmidrule[.08em]{1-5} 
\Block{1-1}
  {
    \diagbox
      {\rule[-2mm]{0pt}{2mm}$S_2$}
      {\rule{0pt}{4mm}$S_1$ }
  }  & $G_1$   & $G_2$   & $\cdots$ & $G_k$  \\
\cmidrule(lr){1-1} \cmidrule(rl){2-5}
$S_{2_1}$              & $\star_{11}$     &   $\star_{12}$  &     $\cdots$    &   $\star_{1k}$ &       \\
$S_{2_2}$              & $\star_{21}$     &   $\star_{22}$  &     $\cdots$    &   $\star_{2k}$ &                   \\
$\vdots$               & $\vdots$ & $\vdots$ & $\ddots $  & $\vdots$ &                \\
$S_{2_\ell}$              & $\star_{\ell1}$     &   $\star_{\ell2}$  &     $\cdots$    & $\star_{\ell k}$                \\   
\cmidrule(lr){1-1} \cmidrule(lr){2-5}
Social coverage $(\cdot100\%)$ & $\frac{\#{\left\{\star_{i1}=\checkmark\right\}}_{i=1}^\ell}{\ell}$ & $\frac{\#{\left\{\star_{i2}=\checkmark\right\}}_{i=1}^\ell}{\ell}$ & $\cdots$ & $\frac{\#{\left\{\star_{ik}=\checkmark\right\}}_{i=1}^\ell}{\ell}$ & \Frame $\frac{\#{\left\{\star_{i j}=\checkmark\right\}}_{i,j=1}^{\ell,k}}{\ell \cdot k}$ \\
\cmidrule[.08em]{1-5} 
\end{NiceTabular}
\end{center}
\end{enumerate}

\paragraph{Total overview of the comparison table} Here, we combine the social and the epidemiological overview of the comparison table.
\begin{enumerate}
    \setcounter{enumi}{5}
    \item Merging of the social and epidemiological overviews.
    \begin{center}
        \begin{NiceTabular}[t]{cccccc}
            \toprule
            \diagbox{$S_2$}{$S_1$} & $G_1$   & $G_2$   & $\cdots$ & $G_k$ & \makecell{Epidemiological \\coverage ($\cdot100\%$)}  \\
            \cmidrule(lr){1-1} \cmidrule(rl){2-5} \cmidrule(rl){6-6} 
            $S_{2_1}$              & $\star_{11}$     &   $\star_{12}$  &     $\cdots$    &   $\star_{1k}$ & $\frac{\#{\left\{\star_{1 j}=\checkmark\right\}}_{j=1}^k}{k}$ \\
            $S_{2_2}$              & $\star_{21}$     &   $\star_{22}$  &     $\cdots$    &   $\star_{2k}$ & $\frac{\#{\left\{\star_{2 j}=\checkmark\right\}}_{j=1}^k}{k}$ \\
            $\vdots$               & $\vdots$ & $\vdots$ & $\ddots $  & $\vdots$ & $\vdots$ \\
            $S_{2_\ell}$              & $\star_{\ell1}$     &   $\star_{\ell2}$  &     $\cdots$    & $\star_{\ell k}$ & $\frac{\#{\left\{\star_{\ell j}=\checkmark\right\}}_{j=1}^k}{k}$ \\   
            \cmidrule(lr){1-1} \cmidrule(lr){2-5}
            Social coverage $(\cdot100\%)$ & $\frac{\#{\left\{\star_{i 1}=\checkmark\right\}}_{i=1}^\ell}{\ell}$ & $\frac{\#{\left\{\star_{i 2}=\checkmark\right\}}_{i=1}^\ell}{\ell}$ & $\cdots$ & $\frac{\#{\left\{\star_{i k}=\checkmark\right\}}_{i=1}^\ell}{\ell}$ & \Frame $\frac{\#{\left\{\star_{i j}=\checkmark\right\}}_{i,j=1}^{\ell,k}}{\ell \cdot k}$ \\
            \cmidrule[.08em]{1-5} 
        \end{NiceTabular}
    \end{center}
\end{enumerate}

\section{Horizontal lockdowns versus age-based interventions} \label{sec:VS}

In this section, we investigate whether age-based interventions can offer a replacement to horizontal lockdowns for the case of SARS-CoV-2, following the framework presented in \hyperref[sec:versus]{\S \ref*{sec:versus}}, and using the model studied in \cite{bgt2023SVEAIRmodel}, which is presented in \hyperref[sec:model]{Appendix \ref*{sec:model}}. In \hyperref[sec:choiceOfScale]{\S \ref*{sec:choiceOfScale}}, we categorize the parameters into $\pmb{\delta}, \pmb{\beta}$ and $\pmb{\gamma}$, and pick our choice of strategic scales, $\pmb{r}$; both $\left(\pmb{\beta},\pmb{\gamma}\right)$ and $\pmb{r}$ serve for the definition of the strategies under investigation. Additionally, we distribute the total population of $\mathscr{P}$ \eqref{SVEIAR-age-scl} into five cohorts, based on age, $\theta$, of each individual. In \hyperref[sec:4Strategies]{\S \ref*{sec:4Strategies}}, we define the graded strategy of horizontal lockdowns and the strategy of age-based restrictions. Finally, in \hyperref[sec:4Simulation]{\S \ref*{sec:4Simulation}}, we compare the aforementioned strategies.


\subsection{Choice of general strategy} \label{sec:choiceOfScale}

The independent variables that appear in $\mathscr{P}$ \eqref{SVEIAR-age-scl} are $t$ and $\theta$, thus $$\pmb{x}=\left(t,\theta\right).$$ In order to define the strategies under investigation, we need to categorize the parameters appeared in $\mathscr{P}$ \eqref{SVEIAR-age-scl} into $(\pmb{\beta}_0, \pmb{\gamma}_0)$ and $\pmb{\delta}$, and consequently choose an appropriate strategic scale as discussed in \hyperref[strategies]{\S \ref*{strategies}}.

The parameters which affect the strategies under investigation are $\beta_A, \beta_I$ and $\gamma_I$. Therefore, we have that
\begin{equation*}
    \left(\pmb{\beta}_0, \gamma_0\right) = \left(\beta_A, \beta_I, \gamma_I\right) \;,
\end{equation*}
whereas
\begin{equation*}
    \pmb{\delta} = \left(\mu, p, \epsilon, \zeta, k, q, \xi, \chi, \gamma_A\right) \;,
\end{equation*}
with the parameter values being as in \hyperref[tab:paramValues]{Table \ref*{tab:paramValues}}.

\begin{table}[H]
\centering
\begin{tabular}{@{}cccc@{}}
\toprule
Parameters & Value                   & Units                         & Source \\ \midrule
$N_0$      & $80 \cdot 10^6$           & individuals                   & Estimated from \cite{ourworldindata}       \\
$\mu$      & $4.38356 \cdot 10^{-5}$ & day $^{-1}$                   & Estimated from \cite{ourworldindata}       \\
$\beta_A$  &  \hyperref[fig:4]{Figure \ref*{fig:4}}                      & individual$^{-1}\, \cdot \, $day $^{-1}$ & Estimated from \cite{DelValle_Hyman_Hethcote_Eubank_2007}       \\
$\beta_I$  &  \hyperref[fig:4]{Figure \ref*{fig:4}}                       & individual$^{-1}\, \cdot \, $day$^{-1}$ & Estimated from \cite{DelValle_Hyman_Hethcote_Eubank_2007}       \\
$p$        & $10^{-3}$               & day$^{-1}$                    & Estimated from \cite{ourworldindata}      \\
$\epsilon$ & 0.7                     & -                             & Estimated from \cite{grant2022impact}      \\
$\zeta$    & $\frac{1}{14}$                    & day$^{-1}$                    & Estimated from \cite{chau2022immunogenicity}      \\
$k$        &  \hyperref[eq:valuek]{Equation \ref*{eq:valuek}}                        & day$^{-1}$                    &  Estimated from \cite{kang2022transmission,wu2022incubation}      \\
$q$        &   \hyperref[fig:5]{Figure \ref*{fig:5}}                       & -                             &  Estimated from \cite{sah2021asymptomatic}      \\
$\xi$      & 0.5                     & -                             & Estimated from \cite{he2021proportion,buitrago2022occurrence}       \\
$\chi$     & \hyperref[eq:valuechi]{Equation \ref*{eq:valuechi}}                         & day$^{-1}$                    & Estimated from \cite{he2021proportion,buitrago2022occurrence}       \\
$\gamma_A$ & $\frac{1}{8}$                    & day$^{-1}$                    & Estimated from \cite{byrne2020inferred}      \\
$\gamma_I$ & $\frac{1}{14}$                    & day$^{-1}$                    & Estimated from \cite{byrne2020inferred}      \\ \bottomrule
\end{tabular}
\caption{A list of parameters of $\mathscr{M}$, along with their value, units, and value source. For their derivation, see \hyperref[sec:paramestim]{Appendix \ref*{sec:paramestim}}.}
\label{tab:paramValues}
\end{table}

We note that regardless of the seemingly important role of asymptomaticity (see \hyperref[sec:model]{Appendix \ref*{sec:model}}) for the spread of the disease, the performance of the means of detection, such as the rapid antigen tests (Ag-RDTs), for the case of asymptomatic infectious individuals still remains ambiguous \citep{CDCP2020,pollock2020asymptomatic,sage56,soni2023performance}. In the light of the above we prefer not to incorporate such means to our general strategy, hence we exclude $\gamma_A$ from $(\pmb{\beta}_0, \pmb{\gamma}_0)$. Moreover, we note that \eqref{uniqcond} holds. 

We are now ready to construct our choice of general strategy along with its strategic scales, following the next steps. 
\begin{enumerate}
\item We consider an interval $\Lambda\subseteq\mathbb{R}_0^+$ such that $0\in \Lambda$, to be the average lifespan of an individual of the population under study, hence $\theta\in\Lambda$. Of course, $\sup{\Lambda}<\infty$.
\item We discretize $\Lambda$ by considering a respective partition $\Delta_\Lambda\coloneqq{\left\{\delta_j\right\}}_{j=0}^n$, for a fixed $n\in\mathbb{N}$, i.e., $$0=\delta_0<\delta_1<\dots<\delta_n=\sup{\Lambda}$$ and we define the subintervals $$\Lambda_j\coloneqq\left[\delta_{j-1},\delta_j\right),\text{ }\forall j\in\left\{1,\dots,n\right\}.$$
\item We set $$\Lambda_W\coloneqq\bigcup\limits_{j \in W}\Lambda_j,\text{ }\forall W\in\mathcal{P}{\left(\left\{1,\dots,n\right\}\right)}\;,$$ where $\mathcal{P}$ stands for the power set, as well as we define $$\begin{aligned}
\rho_W{\left(\,\cdot\,;a\right)}\colon\,\Lambda&\to\left[0,1\right]\\
\theta&\mapsto\rho_W{\left(\theta;a\right)}\coloneqq \begin{cases}
        1   &  \theta \notin \Lambda_W \\
        a   &  \theta \in    \Lambda_W \;,
       \end{cases}
\end{aligned}\text{ }\text{ }\forall\left(W,a\right)\in\mathcal{P}{\left(\left\{1,\dots,n\right\}\right)}\times\left[0,1\right)$$ and 
$$\begin{aligned}
g_W{\left(\,\cdot\,;b\right)}\colon\,\Lambda&\to\left[1,\infty\right]\\
\theta&\mapsto g_W{\left(\theta;b\right)}\coloneqq \begin{cases}
        1   &  \theta \notin \Lambda_W \\
        \frac{1}{b}   &  \theta \in    \Lambda_W \;,
       \end{cases}
\end{aligned}\text{ }\text{ }\forall\left(W,b\right)\in\mathcal{P}{\left(\left\{1,\dots,n\right\}\right)}\times\left[0,1\right),$$ where the choice of $\left(a,b\right)\in{\left[0,1\right)}^2$ is left to be explained. 

We note that in the extreme cases of $W\in\left\{\varnothing,\left\{1,\dots,n\right\}\right\}$ we have $\Lambda_{\varnothing}=\varnothing$ and $\Lambda_{\left\{1,\dots,n\right\}} = \Lambda$, as well as 
\begin{equation*}
    \left(\rho_{\varnothing}(\,\cdot\,; a),g_{\varnothing}(\,\cdot\,; b)\right) = \left(1,1\right)\;,\text{ }\forall\left(a,b\right)\in{\left[0,1\right)}^2 
\end{equation*}
and 
\begin{equation*}
    \left(\rho_{\left\{1,\dots,n\right\}}(\,\cdot\,; a),g_{\left\{1,\dots,n\right\}}(\,\cdot\,; b)\right) = \left(a,\frac{1}{b}\right)\;,\text{ }\forall\left(a,b\right)\in{\left[0,1\right)}^2 . 
\end{equation*}
Hence, the above functions are independent of $\theta$ iff $W\in\left\{\varnothing,\left\{1,\dots,n\right\}\right\}$, as well as they are equal to $1$ iff $W=\varnothing$. 
\item We define the strategic scales to be 
\begin{multline*}
\pmb{r}_{W_{\pmb{\beta}},W_\gamma}{(\,\cdot\,; a,b)} \coloneqq \left( \rho_{W_{\pmb{\beta}}}{(\,\cdot\,;a)}, \rho_{W_{\pmb{\beta}}}{(\,\cdot\,;a)}, g_{W_\gamma}{(\,\cdot\,; b)} \right) \;,\\
\forall\left(W_{\pmb{\beta}},W_\gamma,a,b\right)\in{\left(\mathcal{P}{\left(\left\{1,\dots,n\right\}\right)}\right)}^2\times{\left[0,1\right)}^2 .
\end{multline*}
\item The general strategy of interest $S$ has the form $$S\coloneqq\left\{(\pmb{\beta}, \gamma)=\pmb{r}_{W_{\pmb{\beta}},W_\gamma} (\,\cdot\,; a,b) \odot (\pmb{\beta}_0, \gamma_0)\,\bigg|\,\left(W_{\pmb{\beta}},W_\gamma,a,b\right)\in{\left(\mathcal{P}{\left(\left\{1,\dots,n\right\}\right)}\right)}^2\times{\left[0,1\right)}^2 \right\}.$$
\end{enumerate}

Regarding $a\in\left[0,1\right]$, in the light of \eqref{beta-def}, the effect of every $\rho_{W_{\pmb{\beta}}}(\,\cdot\,,a)$ to $\pmb{\beta}$ can be interpreted as having the average number of close contacts of an (asymptomatic or symptomatic) infectious individual belonging to $\Lambda_{W_{\pmb{\beta}}}$ reduced by $1-a$, i.e. $$\left.\left(c_A,c_I\right)\right|_{\Lambda_{W_{\pmb{\beta}}}}\mapsto a\cdot \left.\left(c_A,c_I\right)\right|_{\Lambda_{W_{\pmb{\beta}}}}.$$ 

Regarding $b\in\left(0,1\right]$, in the light of \eqref{gamma-def}, the effect of every $g_{W_{\gamma}}{(\,\cdot\,,b)}$ to $\gamma$ can be interpreted as having the average period an individual belonging to $\Lambda_{W_{\pmb{\beta}}}$ spends on compartment $I$ before moving into compartment $R$ reduced by  $1-b$, i.e. $$\left.P_{I\to R}\right|_{\Lambda_{W_{\pmb{\beta}}}}\mapsto \left.b\cdot P_{I\to R}\right|_{\Lambda_{W_{\pmb{\beta}}}}.$$ 

\subsection{Choice of distribution of the population into age cohorts}

We now specify the distribution of the whole population into cohorts with respect to the age of each individual, hence with respect to its occupational and social activity. 

\begin{figure}[h]
    \centering
    \includegraphics[width=\textwidth]{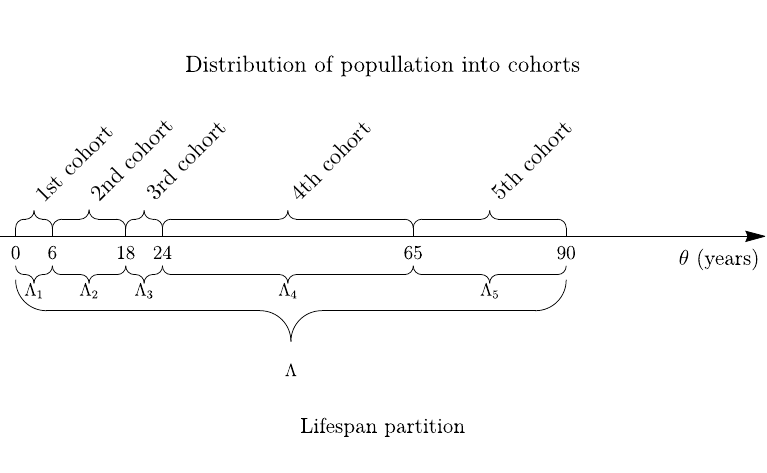}
    \caption{The partition of the non-scaled lifespan and the respective distribution of the whole population into cohorts. The partition was made by taking into account the social profile connecting individuals in each cohort, such as going to school, working, or being pensioners.}
    \label{fig:6}
\end{figure}

We divide the population into five (5) cohorts, thus $n=5$, as seen in \hyperref[fig:6]{Figure \ref*{fig:6}}, where $\Lambda=\left[0,90\right)\;\text{years}$ and $\Delta_\Lambda={\left\{0,6,18,24,65,90\right\}}\;\text{years}$ (both non-scaled). The 1st cohort is made of toddlers and preschoolers, the 2nd is made of school students, the 3rd is primarily made of university students, the 4th is primarily made of the working class, and the 5th is primarily made of pensioners.

As an example of the strategic scale within the context of the cohorts presented in \hyperref[fig:6]{Figure \ref*{fig:6}}, the reduction of the number of contacts of the 1st and 3rd cohort by 80\% and not performing tests on any cohort is modeled by the strategic scale $$\pmb{r}_{\{1,3\},\varnothing}{\left(\,\cdot\,;\frac{1}{5},\,\cdot\,\right)}=\left(\rho_{\{1,3\}}{\left(\,\cdot\,;\frac{1}{5}\right)},g_{\varnothing}{\left(\,\cdot\,; \,\cdot\,\right)}\right),$$ where
\begin{equation*}
    \rho_{\{1,3\}}\left(\theta;\frac{1}{5}\right) = \begin{cases}
    \frac{1}{5},&\text{ if }\theta\in \Lambda_1 \\
    1,&\text{ if }\theta\in \Lambda_2 \\
    \frac{1}{5},&\text{ if }\theta\in \Lambda_3 \\
    1,&\text{ if }\theta\in \Lambda_4 \\
    1,&\text{ if }\theta\in \Lambda_5\; ,
    \end{cases}
\end{equation*}
and
\begin{equation*}
    g_{\varnothing}\left(\theta; \,\cdot\,\right) = 1, \; \forall \theta \in \Lambda \;. 
\end{equation*}
\subsection{Choice of substrategies of general strategy} \label{sec:4Strategies}

In this section, we define the two strategies under investigation. Furthermore, we utilize the strategic scale introduced in \hyperref[sec:choiceOfScale]{\S \ref*{sec:choiceOfScale}} to model each strategy.


\subsubsection{Horizontal lockdowns substrategy} \label{sec:totalNPI}

It is straightforward to check that the largest horizontal with respect to age substrategy of $S$ is 
\begin{multline*}
\bigg\{(\pmb{\beta}, \gamma)=\pmb{r}_{W_{\pmb{\beta}},W_{\gamma}}{\left(\,\cdot\,; a,b\right)} \odot (\pmb{\beta}_0, \gamma_0)\,\bigg|\,\left(W_{\pmb{\beta}},W_{\gamma},a,b\right)\in{\left(\mathcal{P}{\left(\left\{1,\dots,5\right\}\right)}\right)}^2\times{\left[0,1\right)}^2\text{ such that }\\
\text{such that }\left(W_{\pmb{\beta}},W_{\gamma}\right)\in{\left\{\varnothing,\left\{1,\dots,5\right\}\right\}}^2\bigg\},
\end{multline*}
which implies that every substrategy of the above strategy is horizontal with respect to age. 

Thus, a choice of horizontal lockdowns substrategy can be made by considering $S_1\subseteq S$ as $$S_1=\left\{(\pmb{\beta}, \gamma)=\pmb{r}_{\left\{1,\dots,5\right\},\varnothing}{\left(\,\cdot\,; a,\,\cdot\,\right)} \odot (\pmb{\beta}_0, \gamma_0)\,\bigg|\,a\in\left[0,1\right)\right\}.$$ The intensity (that is the amount of contact reduction for every individual) can be varied but uniformly, that is, $$a\in\left[0,1\right)\;\text{ and }\;W_{\pmb{\beta}}=\left\{1,\dots,n\right\},$$ respectively, in order to capture different scenarios. We also assume that no tests are performed in any of the five cohorts, i.e. $$W_\gamma=\varnothing.$$ 

Now, $S_1$ is gradable, since from \eqref{R0-definition} we get that $$\mathcal{R}_0{\left(\pmb{\beta},\gamma,\pmb{\delta}\right)}=a\cdot\mathcal{R}_0{\left(\pmb{\beta}_0,\gamma_0,\pmb{\delta}\right)},\forall\left(\pmb{\beta},\gamma\right)\in S_1\;.$$ In particular, $\mathcal{R}_0$ is strictly increasing with respect to $a$, as it is depicted in \hyperref[tab:R0]{Table \ref*{tab:R0}}.  

\begin{table}[ht]
    \centering
    \begin{tabular}{ccc}
        \toprule
        $a$  &  Contact reduction & $\mathcal{R}_0$ \\
        \midrule
        $0$            & 100\% & 0\\
        $0.1$ & 90\%  & 0.285\\    
        $0.2$ & 80\%  & 0.571\\
        $0.3$ & 70\%  & 0.856\\    
        $0.4$ & 60\%  & 1.141\\
        $0.5$ & 50\%  & 1.427\\
        $0.6$ & 40\%  & 1.712\\
        $0.7$ & 30\%  & 1.998 \\  
        $0.8$ & 20\%  & 2.283\\
        $0.9$& 10\%  & 2.569 \\        
        \bottomrule
    \end{tabular}
    \caption{The value of $\mathcal{R}_0$ decreases linearly as the intensity of the stay-at-home restrictions increases, i.e. as $a$ decreases.}
    \label{tab:R0}
\end{table}

To get a better understanding of how $a\in\left[0,1\right)$ translates into the real life intensity of a stay-at-home restriction policy, we firstly notice that when $W_{\pmb{\beta}}=\varnothing$ (or else $a\to 1^-$), we have that no stay-at-home restrictions are in effect. In that case, our model predicts an $\mathcal{R}_0$ of 2.854 (or else $\mathcal{R}_0\to 2.854^-$), which is in line with various systematic reviews found in scientific literature, such as 2.87 (95\% CI: 2.39-–3.44) in \cite{billah2020reproductive} and 2.69 (95\% CI: 2.40--2.98) in \cite{ahammed2021estimation}, which solidifies the validity of our model in predicting the $\mathcal{R}_0$ of SARS-CoV-2 pandemic. Furthermore, we see that the tight lockdown Italy enforced in early 2020 resulted in an 82\% reduction in mobility \citep{vinceti2022substantial}, which would correspond to $a$ being approximately equal to 0.2. During the same time frame in Germany, the authors of \cite{schlosser2020covid}, report about a 50\% drop in the average number of contacts, which corresponds to $a = 0.5$. Finally, in \cite{zhou2020effects} the authors show that even a 20\% reduction in mobility proved a good way of delaying the spread of the infection, which would correspond to $a$ being approximately equal to 0.8.

Based on the aforementioned cases, we construct three different scenarios based on the intensity of the mobility restrictions:
\begin{itemize}
    \item the low intensity scenario, $\mathcal{L}$, where the average number of contacts is reduced by 20\% and ${\mathcal{R}_0}_{\mathcal{L}} = 2.283$,
    \item the medium intensity scenario, $\mathcal{M}$, where the average number of contacts is reduced by 50\% and ${\mathcal{R}_0}_{\mathcal{M}} = 1.427$ and 
    \item the high intensity scenario, $\mathcal{H}$, where the average number of contacts is reduced by 80\% and ${\mathcal{R}_0}_{\mathcal{H}} = 0.571$.
\end{itemize}

The above scenarios constitute a gradation $\pmb{G}$ of $S_1$ with $$G_1={\mathcal{R}_0}_{\mathcal{H}},\text{ }G_2={\mathcal{R}_0}_{\mathcal{M}}\text{ and }G_3={\mathcal{R}_0}_{\mathcal{L}}.$$ Such a gradation is summarized in  \hyperref[tab:totalLockdownSen]{Table \ref*{tab:totalLockdownSen}}.

\begin{table}[h]
    \centering
    \begin{tabular}{cccc}
        \toprule
         Gradation & Intensity level                & Contact reduction & $\mathcal{R}_0$ \\ \midrule
         $G_1$ & \textcolor{red}{High ($\mathcal{H}$)}      & 80\% & 0.571 \\
         $G_2$ & \textcolor{orange}{Medium ($\mathcal{M}$)}  & 50\% & 1.427 \\
         $G_3$ & \textcolor{magenta}{Low ($\mathcal{L}$)}    & 20\% & 2.283 \\ \bottomrule
    \end{tabular}
    \caption{Summary of the three horizontal lockdowns' intensity scenarios, Low ($\mathcal{L}$), Medium ($\mathcal{M}$) and High ($\mathcal{H}$), which constitute a gradation of $S_1$.}
    \label{tab:totalLockdownSen}
\end{table}

\subsubsection{Aged-based substrategies}\label{sec:agedbasedNPI} 

The largest aged-based substrategy of $S$, $S_2\subseteq S$, is 
\begin{multline*}
S_2=\bigg\{(\pmb{\beta}, \gamma)=\pmb{r}_{W_{\pmb{\beta}},W_{\gamma}}{\left(\,\cdot\,; a,b\right)} \odot (\pmb{\beta}_0, \gamma_0)\,\bigg|\,\left(W_{\pmb{\beta}},W_{\gamma},a,b\right)\in{\left(\mathcal{P}{\left(\left\{1,\dots,5\right\}\right)}\right)}^2\times{\left[0,1\right)}^2\text{ such that }\\
\text{such that }\left(W_{\pmb{\beta}},W_{\gamma}\right)\nin{\left\{\varnothing,\left\{1,\dots,5\right\}\right\}}^2\bigg\},
\end{multline*}
which implies that every substrategy of $S_2$ is aged-based hence it can potentially be compared to the graded $S_1$. For our simulations, we choose the substrategies ${\left\{S_{2_i}\right\}}_{i=1}^{16}$ of $S_2$ summarized in \hyperref[tab:agebasedStrategies]{Table \ref*{tab:agebasedStrategies}} for the comparison to $S_1$.

\begin{table}[H]
\centering
\begin{tabular}{@{}cccc@{}}
\toprule
$i$ & Age-based interventions & $W_{\pmb{\beta}}$ & $W_{\gamma}$\\ \midrule
1   & \begin{tabular}[c]{@{}c@{}}Contact reduction: 1st, 2nd, 3rd cohorts \\  Testing: 4th, 5th cohorts\end{tabular} & $\left\{1,2,3\right\}$ & $\left\{4,5\right\}  $ \\ 
2   & \begin{tabular}[c]{@{}c@{}}Contact reduction: 4th, 5th cohorts \\  Testing: 1st, 2nd, 3rd cohorts\end{tabular} & $\left\{4,5\right\}$   & $\left\{1,2,3\right\}$ \\ 
3   & \begin{tabular}[c]{@{}c@{}}Contact reduction: 1st cohort \\  Testing: 4th, 5th cohorts\end{tabular} & $\left\{1\right\}$   & $\left\{4,5\right\}$  \\               
4   & \begin{tabular}[c]{@{}c@{}}Contact reduction: 4th, 5th cohorts \\  Testing: 1st cohort\end{tabular} & $\left\{4,5\right\}$ & $\left\{1\right\}  $   \\              
5   & \begin{tabular}[c]{@{}c@{}}Contact reduction: 2nd cohort \\  Testing: 4th, 5th cohorts\end{tabular} & $\left\{2\right\}$   & $\left\{4,5\right\}$  \\               
6   & \begin{tabular}[c]{@{}c@{}}Contact reduction: 4th, 5th cohorts \\  Testing: 2nd cohort\end{tabular} & $\left\{4,5\right\}$ & $\left\{2\right\}  $  \\               
7   & \begin{tabular}[c]{@{}c@{}}Contact reduction: 3rd cohort \\  Testing: 4th, 5th cohorts\end{tabular} & $\left\{3\right\}$   & $\left\{4,5\right\}$  \\               
8   & \begin{tabular}[c]{@{}c@{}}Contact reduction: 4th, 5th cohorts \\  Testing: 3rd cohort\end{tabular} & $\left\{4,5\right\}$ & $\left\{3\right\}  $  \\               
9   & \begin{tabular}[c]{@{}c@{}}Contact reduction: 1st cohort \\  Testing: 2nd cohort\end{tabular}       & $\left\{1\right\}$ & $\left\{2\right\}$ \\                    
10  & \begin{tabular}[c]{@{}c@{}}Contact reduction: 2nd cohort \\  Testing: 1st cohort\end{tabular}       & $\left\{2\right\}$ & $\left\{1\right\}$  \\                   
11  & \begin{tabular}[c]{@{}c@{}}Contact reduction: 4th cohort \\  Testing: 5th cohort\end{tabular}       & $\left\{4\right\}$ & $\left\{5\right\}$ \\                    
12  & \begin{tabular}[c]{@{}c@{}}Contact reduction: 5th cohort \\  Testing: 4th cohort\end{tabular}       & $\left\{5\right\}$ & $\left\{4\right\}$ \\                    
13  & \begin{tabular}[c]{@{}c@{}}Contact reduction: 2nd cohort \\  Testing: 4th cohort\end{tabular}       & $\left\{2\right\}$ & $\left\{4\right\}$  \\                   
14  & \begin{tabular}[c]{@{}c@{}}Contact reduction: 4th cohort \\  Testing: 2nd cohort\end{tabular}       & $\left\{4\right\}$ & $\left\{2\right\}$   \\                  
15  & \begin{tabular}[c]{@{}c@{}}Contact reduction: 2nd cohort \\  Testing: 5th cohort\end{tabular}       & $\left\{2\right\}$ & $\left\{5\right\}$ \\                    
16  & \begin{tabular}[c]{@{}c@{}}Contact reduction: 5th cohort \\  Testing: 2nd cohort\end{tabular}       & $\left\{5\right\}$ & $\left\{2\right\}$ \\ \bottomrule
\end{tabular} 

\caption{The sixteen age-based substrategies ${\left\{S_{2_i}\right\}}_{i=1}^{16}$ of $S_2$ which are chosen for the comparison to $S_1$.}
\label{tab:agebasedStrategies}
\end{table}


\subsection{Simulations and results}  \label{sec:4Simulation}

Here we employ the scheme introduced in \hyperref[comparison]{\S \ref*{comparison}} for the comparison between $S_1$ of \hyperref[sec:totalNPI]{\S \ref*{sec:totalNPI}} and ${\left\{S_{2_i}\right\}}_{i=1}^{16}$ of \hyperref[sec:agedbasedNPI]{\S \ref*{sec:agedbasedNPI}}. The simulations were performed using Mathematica 13.1 \cite{Mathematica}. In the end of this section, we summarize its results with the comparison table. 

\begin{figure}[ht]
     \centering
     \begin{subfigure}[b]{0.48\textwidth}
         \centering
         \includegraphics[width=\textwidth]{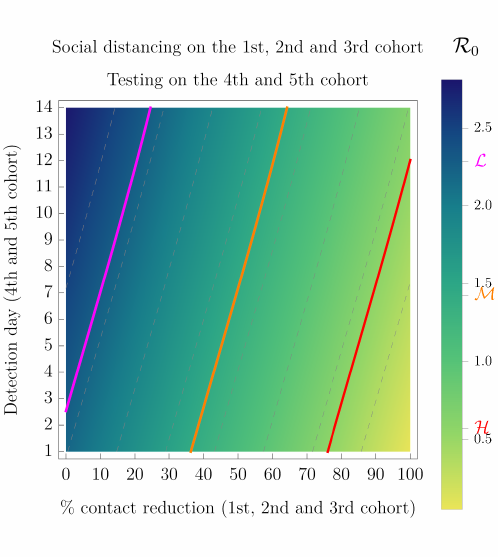}
         \caption{$\mathcal{R}_0$ when social distancing is enforced on the 1st, 2nd and 3rd cohort and testing is enforced on the 5th and 6th cohort ($W_{\pmb{\beta}} = \left\{1,2,3\right\}$ and $W_{\gamma} = \left\{4,5\right\}$).}
         \label{fig:7a}
     \end{subfigure}
     \hfill
     \begin{subfigure}[b]{0.48\textwidth}
         \centering
         \includegraphics[width=\textwidth]{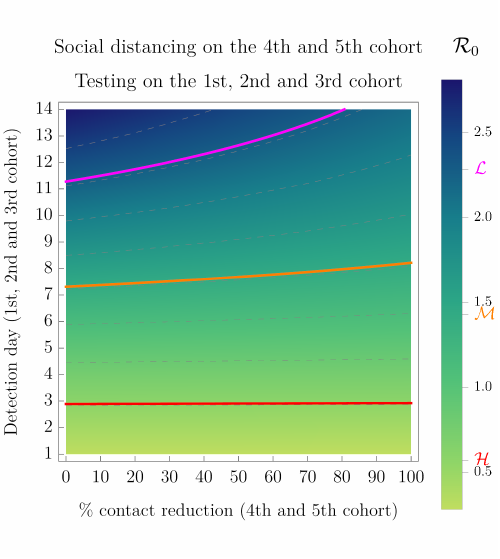}
         \caption{$\mathcal{R}_0$ when social distancing is enforced on the 5th and 6th cohort and testing is enforced on the 1st, 2nd and 3rd cohort ($W_{\pmb{\beta}} = \left\{5, 6\right\}$ and $W_{\gamma} = \left\{1, 2, 3\right\}$).}
         \label{fig:7b}
     \end{subfigure}
        \caption{Two density plots of the grouping of the three younger cohorts and the two older cohorts together. In both cases, all three of our horizontal lockdown scenarios $\mathcal{L}, \mathcal{M}$ and $\mathcal{H}$, can be replaced by enforcing a wide range of intensity level restrictions to the different cohort groupings.}
        \label{fig:7}
\end{figure}


\subsubsection{Social overview of the results, \texorpdfstring{$S_2$}{S2} versus \texorpdfstring{$S_1$}{S1}}
\label{Social overview of the results} 

Throughout our simulations, we let the $\left(a,b\right)$ of each strategic scale to take values in the 2D interval $[0,1)^2$ and illustrate the  results in density plots, where in the $x$-axis and $y$-axis we have $a\cdot 100\%$ and $b\cdot P_{I\to R}=\frac{b}{\gamma_I}=b\cdot 14$ days, respectively. 

\paragraph{$S_{2_{1,2}}$ versus $S_1$} We begin by examining whether restrictions on the younger or the older cohorts play a more important role in reducing $\mathcal{R}_0$. In \hyperref[fig:7a]{Figure \ref*{fig:7a}}, we see that in order to achieve the same $\mathcal{R}_0$ as the scenario $\mathcal{H}$, the contact reduction of the first three cohorts needs to be at least 75\% and the individuals of the last two cohorts need to be detected and removed at least before the twelfth day. Additionally, since the absolute value of the gradient of the contour lines is high, the younger cohorts influence the dynamics of $\mathcal{R}_0$ more when compared to the older cohorts. In \hyperref[fig:7b]{Figure \ref*{fig:7b}}, we see that the scenario $\mathcal{H}$, can be replaced by finding and removing from the community the people belonging in the first three cohorts at around the third day from symptom onset, whereas the contact reduction of the older age cohorts is almost irrelevant. Furthermore, since the gradient of the contour lines is almost zero, the younger cohorts play a far greater role in reducing $\mathcal{R}_0$ when compared to the older cohorts, especially the more austere the restrictions are. Overall, \hyperref[fig:7]{Figure \ref*{fig:7}} shows us that the younger cohorts are more influential in the dynamics of $\mathcal{R}_0$, both when they are faced with social distancing restrictions, and with mandatory testing.

It is now clear that the younger cohorts play a far more important role in the dynamics of $\mathcal{R}_0$. We subsequently examine whether similar results as those presented in \hyperref[fig:7]{Figure \ref*{fig:7}} can be achieved, by restricting just one of the three younger cohorts instead of all three of them together.

\begin{figure}[!h]
     \centering
     \begin{subfigure}[b]{0.48\textwidth}
            \centering
            \includegraphics[width=\textwidth]{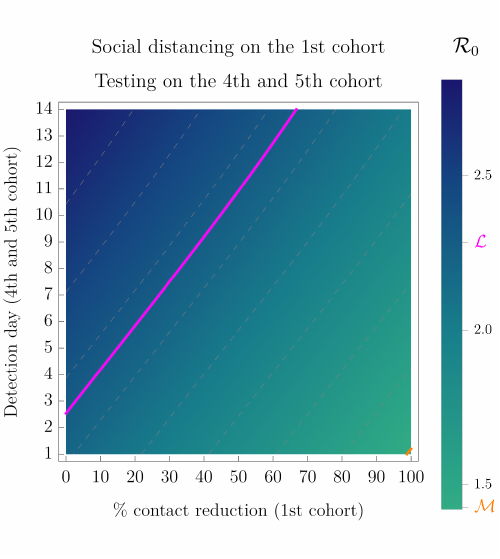}
            \caption{$\mathcal{R}_0$ when social distancing is enforced on the 1st cohort and testing is enforced on the 5th and 6th cohort ($W_{\pmb{\beta}} = \left\{1\right\}$ and $W_{\gamma} = \left\{4,5\right\}$).}
            \label{fig:8a}
     \end{subfigure}
     \hfill
     \begin{subfigure}[b]{0.48\textwidth}
            \centering
            \includegraphics[width=\textwidth]{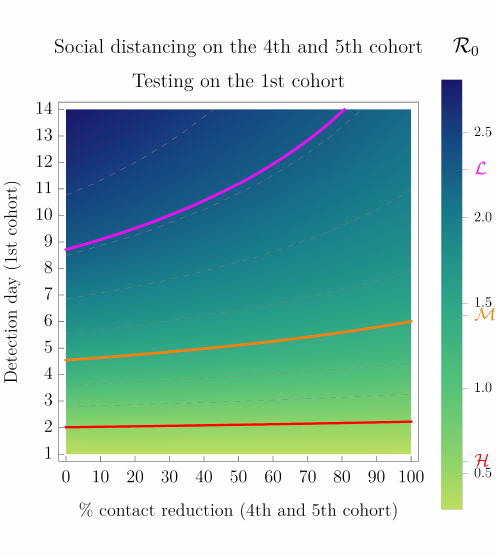}
            \caption{$\mathcal{R}_0$ when social distancing is enforced on the 5th and 6th cohort and testing is enforced on the 1st cohort ($W_{\pmb{\beta}} = \left\{5, 6\right\}$ and $W_{\gamma} = \left\{1\right\}$).}
            \label{fig:8b}
     \end{subfigure}
        \caption{Two density plots of the 1st cohort and the grouping of the two older cohorts together. When social distancing is enforced on the 1st cohort scenario $\mathcal{H}$ can only be achieved when the most austere restriction are enforced. On the other hand, when the symptomatic individuals of the 1st cohort are the ones getting tested all three of our horizontal lockdown scenarios $\mathcal{L}, \mathcal{M}$ and $\mathcal{H}$, can be replaced by enforcing a wide range of intensity level restrictions to the 1st cohort and the grouping of the 4th and 5th cohort. The detection-and-removal day of asymptomatic individuals needs to be one day faster when compared to the simulation illustrated in \hyperref[fig:8b]{Figure \ref*{fig:8b}}, for the same results as scenario $\mathcal{H}$ to apply.}
        \label{fig:8}
\end{figure}

\paragraph{$S_{2_{3,4}}$ versus $S_1$} \hyperref[fig:8]{Figure \ref*{fig:8}} illustrates restrictions on the 1st and the 4th -- 5th cohorts. When social distancing on the 1st cohort and testing on the 4th and 5th cohorts are enforced, scenario $\mathcal{M}$ can only be achieved with the strongest possible restrictions on the aforementioned cohorts, as we can see in \hyperref[fig:8a]{Figure \ref*{fig:8a}}. When the restrictions on the cohorts are reversed, \hyperref[fig:8b]{Figure \ref*{fig:8b}} shows that scenario $\mathcal{H}$ can be achieved if the day that the symptomatic infectious individuals are detected and removed from the community is around the second day, with the contact reduction of the 4th and 5th cohort being almost irrelevant just like the case described by \hyperref[fig:7]{Figure \ref*{fig:7}}. There is, however, a one-day difference in the required detection day of asymptomatic individuals between the scenarios presented in \hyperref[fig:7b]{Figure \ref*{fig:7b}} and \hyperref[fig:8b]{Figure \ref*{fig:8b}} for them to have the same effect on $R_0$, as scenario $\mathcal{H}$. In other words, the procedure of detection and removal of asymptomatic individuals from the community needs to be one day faster when only the 1st cohort is getting tested when compared to the grouping of the 1st, 2nd and 3rd cohorts, for them to have the same results on $\mathcal{R}_0$ as scenario $\mathcal{H}$. 

\begin{figure}[!h]
     \centering
     \begin{subfigure}[b]{0.48\textwidth}
            \centering
            \includegraphics[width=\textwidth]{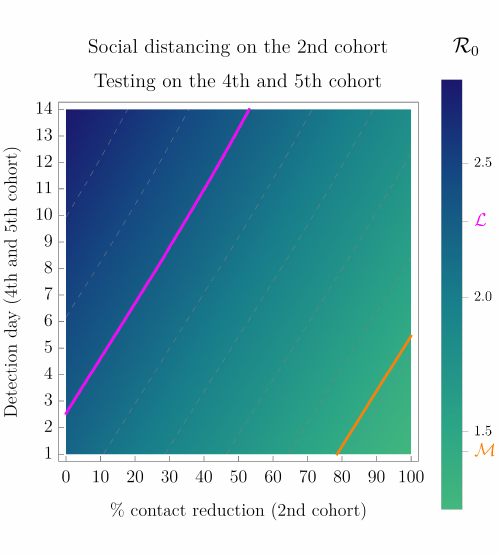}
            \caption{$\mathcal{R}_0$ when social distancing is enforced on the 2nd cohort and testing is enforced on the 5th and 6th cohort ($W_{\pmb{\beta}} = \left\{2\right\}$ and $W_{\gamma} = \left\{4,5\right\}$).}
            \label{fig:9a}
     \end{subfigure}
     \hfill
     \begin{subfigure}[b]{0.48\textwidth}
            \centering
            \includegraphics[width=\textwidth]{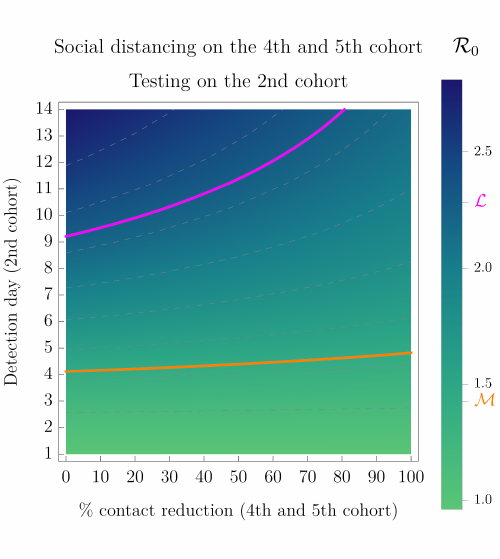}
            \caption{$\mathcal{R}_0$ when social distancing is enforced on the 5th and 6th cohort and testing is enforced on the 2nd cohort ($W_{\pmb{\beta}} = \left\{5, 6\right\}$ and $W_{\gamma} = \left\{2\right\}$).}
            \label{fig:9b}
     \end{subfigure}
        \caption{Two density plots of the 2nd cohort and the grouping of the two older cohorts together. When social distancing is enforced on the 2nd cohort, scenario $\mathcal{H}$ can only be achieved with laxer restriction compared to the respective restrictions on the 1st cohort. Neither of the pictured simulations are able to offer a replacement to scenario $\mathcal{H}$. Much like the simulations of \hyperref[fig:7b]{Figure \ref*{fig:7b}} and \hyperref[fig:8b]{Figure \ref*{fig:8b}}, for the scenario $\mathcal{M}$ to be achieved the testing of the younger cohorts dominates the dynamics of $\mathcal{R}_0$, with the dynamics of the older cohorts being almost irrelevant.}
        \label{fig:9}
\end{figure}

\paragraph{$S_{2_{5,6}}$ versus $S_1$} Next, we examine the importance of the 2nd cohort to the dynamics of $\mathcal{R}_0$, with the results being shown in \hyperref[fig:9]{Figure \ref*{fig:9}}. Contrary to the simulation of \hyperref[fig:8a]{Figure \ref*{fig:8a}}, when social distance is enforced on the 2nd cohort, the results of scenario $\mathcal{M}$ can be achieved with far less strict policies. In particular, as shown in \hyperref[fig:9a]{Figure \ref*{fig:9a}}, for scenario $\mathcal{M}$ to be achieved, the contacts of the 2nd cohort need to be reduced by at least 80\% and the symptomatic individuals of the 4th and 5th cohorts need to be detected and removed from the community at least before around the fifth day. When the 2nd cohort is the one being tested, \hyperref[fig:9b]{Figure \ref*{fig:9b}} shows that scenario $\mathcal{M}$ can be achieved by removing symptomatic individuals from the community at around the fifth day, with the reduction in the average number of contacts of the 4th and 5th cohorts being almost irrelevant, much like the simulations illustrated in \hyperref[fig:7b]{Figure \ref*{fig:7b}} and \hyperref[fig:8b]{Figure \ref*{fig:8b}}. Additionally, none of the simulations of \hyperref[fig:9]{Figure \ref*{fig:9}} can act as a replacement measure to scenario $\mathcal{H}$.

\begin{figure}[!h]
     \centering
     \begin{subfigure}[b]{0.48\textwidth}
            \centering
            \includegraphics[width=\textwidth]{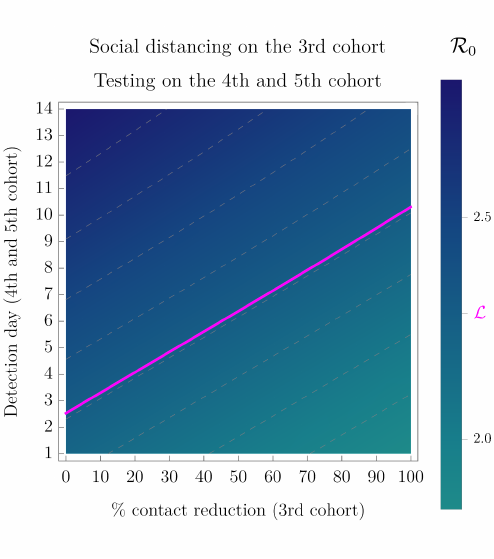}
            \caption{$\mathcal{R}_0$ when social distancing is enforced on the 3rd cohort and testing is enforced on the 5th and 6th cohort ($W_{\pmb{\beta}} = \left\{3\right\}$ and $W_{\gamma} = \left\{4,5\right\}$).}
            \label{fig:10a}
     \end{subfigure}
     \hfill
     \begin{subfigure}[b]{0.48\textwidth}
            \centering
            \includegraphics[width=\textwidth]{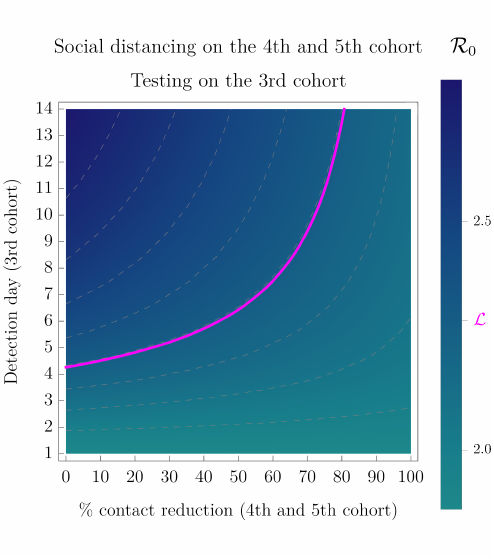}
            \caption{$\mathcal{R}_0$ when social distancing is enforced on the 5th and 6th cohort and testing is enforced on the 3rd cohort ($W_{\pmb{\beta}} = \left\{5, 6\right\}$ and $W_{\gamma} = \left\{3\right\}$).}
            \label{fig:10b}
     \end{subfigure}
        \caption{Two density plots of the 3rd cohort and the grouping of the two older cohorts together. Neither of the simulations is able to offer a replacement to scenario $\mathcal{M}$ and scenario $\mathcal{H}$. The influence of the 3rd cohort to the dynamics of $\mathcal{R_0}$ is far weaker when compared to the influence of the 1st and 2nd cohort, as can be seen from \hyperref[fig:8]{Figure \ref*{fig:8}} and \hyperref[fig:9]{Figure \ref*{fig:9}}.}
        \label{fig:10}
\end{figure}

\paragraph{$S_{2_{7,8}}$ versus $S_1$} Subsequently, we examine the contribution of the 3rd cohort to the dynamics of $\mathcal{R_0}$, with \hyperref[fig:10]{Figure \ref*{fig:10}} illustrating the results. As we can see from \hyperref[fig:10]{Figure \ref*{fig:10}}, the 3rd cohort, in combination with the grouping of the 4th and 5th cohort, seems to influence the reduction of $\mathcal{R}_0$ far less when compared to the younger cohorts. The only horizontal lockdown scenario that can be replaced with this combination of age-based interventions is scenario $\mathcal{L}$. Additionally, even though the 1st and 2nd cohort dominated the dynamics of $\mathcal{R}_0$ when the symptomatic individuals of those cohorts were getting tested, that is not the case with the 3rd cohort, as can be seen from \hyperref[fig:10b]{Figure \ref*{fig:10b}}. The same holds for the case when social distancing is enforced on the 3rd cohort, since the absolute value of the gradient of the contour lines of \hyperref[fig:10a]{Figure \ref*{fig:10a}} is about 2. Hence, out of the three younger cohorts, the 3rd one has the weakest influence on the dynamics of $\mathcal{R}_0$.

\begin{figure}[ht]
     \centering
     \begin{subfigure}[b]{0.48\textwidth}
            \centering
            \includegraphics[width=\textwidth]{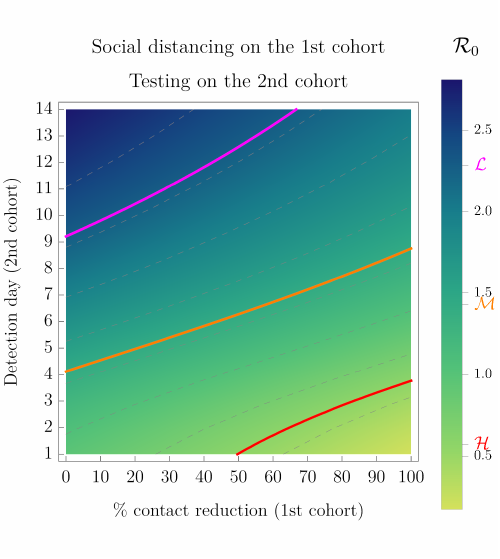}
            \caption{$\mathcal{R}_0$ when social distancing is enforced on the 1st cohort and testing is enforced on the 2nd cohort ($W_{\pmb{\beta}} = \left\{1\right\}$ and $W_{\gamma} = \left\{2\right\}$).}
            \label{fig:11a}
     \end{subfigure}
     \hfill
     \begin{subfigure}[b]{0.48\textwidth}
            \centering
            \includegraphics[width=\textwidth]{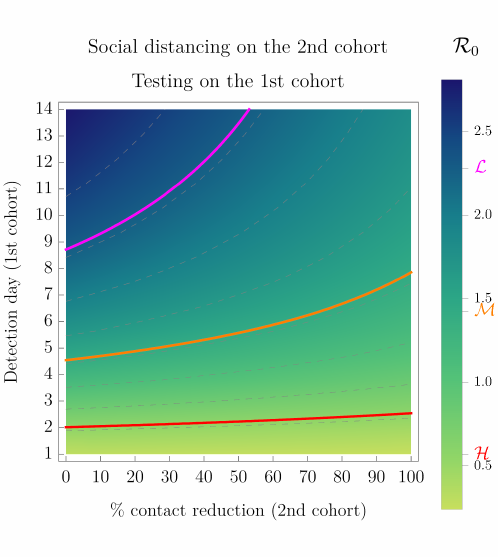}
            \caption{$\mathcal{R}_0$ when social distancing is enforced on the 2nd cohort and testing is enforced on the 1st cohort ($W_{\pmb{\beta}} = \left\{2\right\}$ and $W_{\gamma} = \left\{1\right\}$).}
            \label{fig:11b}
     \end{subfigure}
        \caption{Two density plots of the influence of the 1st cohort and 2nd cohort on the dynamics of $\mathcal{R}_0$. In both cases, all three of our horizontal lockdown scenarios $\mathcal{L}, \mathcal{M}$ and $\mathcal{H}$, can be replaced by enforcing a wide range of intensity level restrictions to the 1st cohort and 2nd cohort. The 1st and 2nd cohorts are the most important cohorts at effecting the dynamics of $\mathcal{R}_0$, since they influence the dynamics of $\mathcal{R}_0$ comparably to the influence of the combination of all of our cohorts, as seen in \hyperref[fig:8]{Figure \ref*{fig:8}}.}
        \label{fig:11}
\end{figure}

\paragraph{$S_{2_{9,10}}$ versus $S_1$} The 1st and 2nd cohort seem to be the two cohorts that influence the dynamics of $\mathcal{R}_0$ the most. Hence, we quantify the results of targeting only the aforementioned cohorts in \hyperref[fig:11]{Figure \ref*{fig:11}}. As we can see from \hyperref[fig:11]{Figure \ref*{fig:11}}, all three horizontal lockdown scenarios $\mathcal{L}, \mathcal{M}$ and $\mathcal{H}$ can be replaced with a combination of measures targeted at the 1st and 2nd cohort. This particular combination of age-based measures has similar dynamics as the scenarios presented in \hyperref[fig:7]{Figure \ref*{fig:7}} which combine all of our cohorts, and \hyperref[fig:8b]{Figure \ref*{fig:8b}} which includes measures regarding three different cohorts. The vital role of the 1st and 2nd cohort is now undeniable. In \hyperref[fig:11a]{Figure \ref*{fig:11a}} we see, that scenario $\mathcal{H}$ can be replaced with the contacts of the 1st cohort being reduced by at least 50\% and the infectious individuals of the 2nd cohort being found and removed from the community at least before the 4th day. When the restrictions are reversed, scenario $\mathcal{H}$ can be replaced when the symptomatic individuals of the 1st cohort are detected and removed from the community at around the second day after symptom onset, as can be seen in \hyperref[fig:11b]{Figure \ref*{fig:11b}}, with minimal contribution from the 2nd cohort.

\begin{figure}[!h]
     \centering
     \begin{subfigure}[b]{0.48\textwidth}
            \centering
            \includegraphics[width=\textwidth]{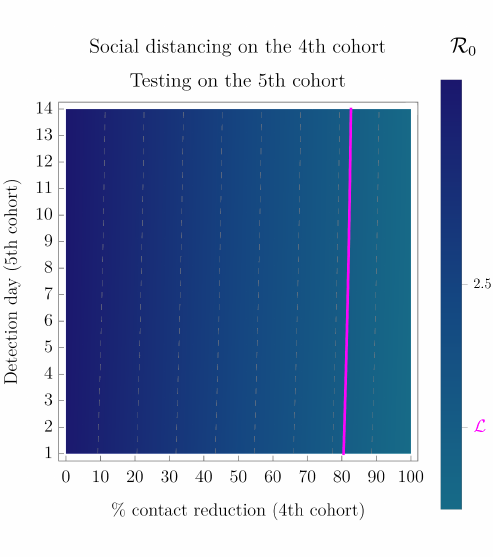}
            \caption{$\mathcal{R}_0$ when social distancing is enforced on the 4th cohort and testing is enforced on the 5th cohort ($W_{\pmb{\beta}} = \left\{4\right\}$ and $W_{\gamma} = \left\{5\right\}$).}
            \label{fig:12a}
     \end{subfigure}
     \hfill
     \begin{subfigure}[b]{0.48\textwidth}
            \centering
            \includegraphics[width=\textwidth]{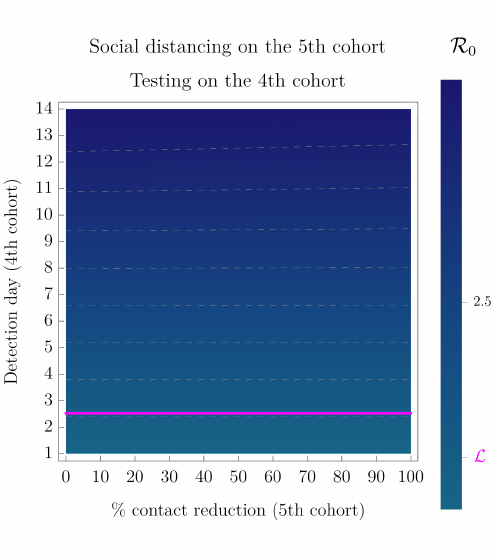}
            \caption{$\mathcal{R}_0$ when social distancing is enforced on the 5th cohort and testing is enforced on the 4th cohort ($W_{\pmb{\beta}} = \left\{5\right\}$ and $W_{\gamma} = \left\{4\right\}$).}
            \label{fig:12b}
     \end{subfigure}
        \caption{Two density plots of the influence of the 4th and 5th cohort on the dynamics of $\mathcal{R}_0$. Neither case was able to offer a replacement to horizontal lockdown scenario $\mathcal{M}$ and scenario $\mathcal{H}$. Restrictions on the combination of the 4th cohort and the 5th cohort result in the poorest reduction in $\mathcal{R}_0$ when compared to the remaining of our simulations. When measures are imposed to the 4th and 5th cohort, the restrictions on the 4th cohort dominate the dynamics of $\mathcal{R}_0$.}
        \label{fig:12}
\end{figure}

\paragraph{$S_{2_{11,12}}$ versus $S_1$} Up until now, we examined the two older cohorts, namely the 4th and 5th cohort, grouping them together as a single cohort. In an attempt to study the result of the interactions of the aforementioned cohorts individually, we present \hyperref[fig:12]{Figure \ref*{fig:12}}. As expected, from the inability of the grouping of the 4th and 5th cohort to dominate the dynamics of our previous simulations, the simulations of \hyperref[fig:12]{Figure \ref*{fig:12}} offer a poor reduction of $\mathcal{R}_0$. Neither in \hyperref[fig:12a]{Figure \ref*{fig:12a}} nor \hyperref[fig:12b]{Figure \ref*{fig:12b}} can horizontal lockdown scenarios $\mathcal{H}$ and $\mathcal{M}$ be replaced by a combination of measures in the 4th and 5th cohort. Only scenario $\mathcal{L}$ can be replaced, and that is with austere restrictions on the 5th cohort. In particular, scenario $\mathcal{L}$ can be achieved either when the reduction of the average amount of contacts of the 4th cohort is 80\% or when the symptomatic individuals of the 4th cohort are detected and removed from the community at around 2.5 days after symptom onset. Finally, there is a clear domination of the 4th cohort in this particular combination of age-based measures, with the measures enforced on the 5th cohort being irrelevant.

\begin{figure}[!h]
     \centering
     \begin{subfigure}[b]{0.43\textwidth}
            \centering
          \includegraphics[width=\textwidth]{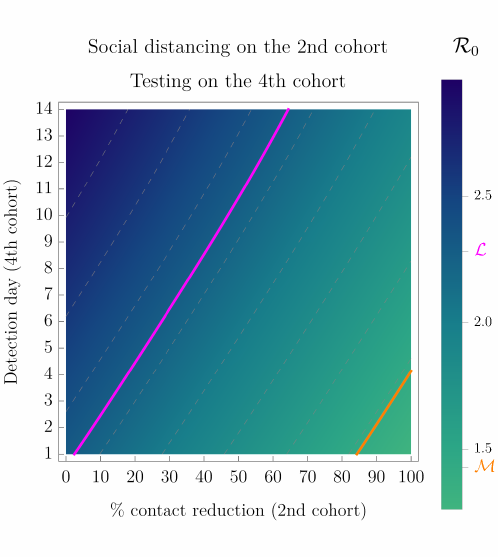}
            \caption{$\mathcal{R}_0$ when social distancing is enforced on the 2nd cohort and testing is enforced on the 4th cohort ($W_{\pmb{\beta}} = \left\{2\right\}$ and $W_{\gamma} = \left\{4\right\}$).}
            \label{fig:13a}
     \end{subfigure}
     \hfill
     \begin{subfigure}[b]{0.43\textwidth}
            \centering
            \includegraphics[width=\textwidth]{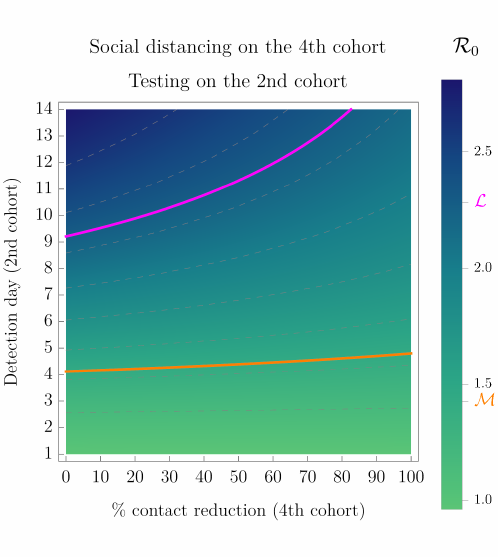}
            \caption{$\mathcal{R}_0$ when social distancing is enforced on the 4th cohort and testing is enforced on the 2nd cohort ($W_{\pmb{\beta}} = \left\{4\right\}$ and $W_{\gamma} = \left\{2\right\}$).}
            \label{fig:13b}
     \end{subfigure}
             \vskip\baselineskip
     \begin{subfigure}[b]{0.43\textwidth}
            \centering
            \includegraphics[width=\textwidth]{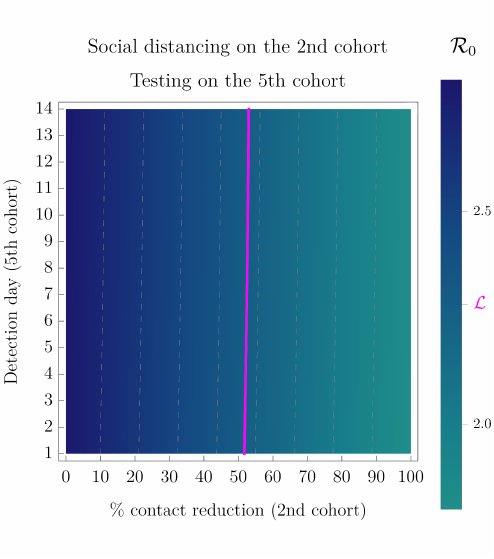}
            \caption{$\mathcal{R}_0$ when social distancing is enforced on the 2nd cohort and testing is enforced on the 5th cohort ($W_{\pmb{\beta}} = \left\{2\right\}$ and $W_{\gamma} = \left\{5\right\}$).}
            \label{fig:13c}
     \end{subfigure}
     \hfill
     \begin{subfigure}[b]{0.43\textwidth}
            \centering
            \includegraphics[width=\textwidth]{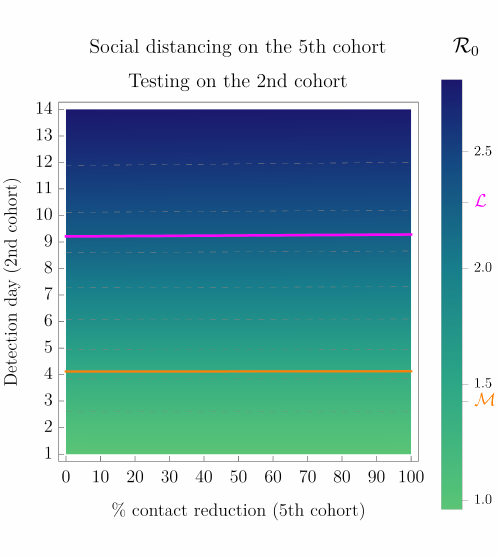}
            \caption{$\mathcal{R}_0$ when social distancing is enforced on the 5th cohort and testing is enforced on the 2nd cohort ($W_{\pmb{\beta}} = \left\{5\right\}$ and $W_{\gamma} = \left\{2\right\}$).}
            \label{fig:13d}
     \end{subfigure}
        \caption{Four density plots of the influence of the interactions of the 2nd and 4th cohort, as well as the 2nd and 5th cohort, on the dynamics of $\mathcal{R}_0$. None of the simulations was able to offer a replacement to horizontal lockdown scenario $\mathcal{H}$. The 5th cohort's contribution to the dynamics of $\mathcal{R}_0$ is insignificant, since its removal from the measure-targeted cohorts, minimally affects the dynamics of $\mathcal{R}_0$, as can be seen when \hyperref[fig:13a]{Figure \ref*{fig:13a}}, \hyperref[fig:13b]{Figure \ref*{fig:13b}} and \hyperref[fig:9]{Figure \ref*{fig:9}} are compared. Additionally, the 2nd cohort completely dominates the dynamics of $\mathcal{R}_0$, when the 5th cohort is included in the simulations, as can be seen from \hyperref[fig:13c]{Figure \ref*{fig:13c}} and \hyperref[fig:13d]{Figure \ref*{fig:13d}}.}
        \label{fig:13}
\end{figure}

\paragraph{$S_{2_{13,14,15,16}}$ versus $S_1$} Lastly, we present the final combination of measures in \hyperref[fig:13]{Figure \ref*{fig:13}}. This final set of restrictions acts as a viable proposal to a real life situation with the economic impact of the measures in mind, since it targets the 2nd cohort, i.e., school students, whose contact reduction, or in other words school closures, would minimally affect the economy. Additionally, \hyperref[fig:13]{Figure \ref*{fig:13}}, allows us to examine the difference between the grouping of the two older cohorts and their individual contribution to $\mathcal{R}_0$, in combination to another, younger, cohort. As can be seen in \hyperref[fig:13a]{Figure \ref*{fig:13a}}, for horizontal lockdown scenario $\mathcal{M}$ to be replaced, the contact reduction of the 2nd cohort needs to be at least 85\% and the infectious individuals of the 4th cohort need to be found and removed from the community at least before the fourth day after symptom onset. Compared to \hyperref[fig:9a]{Figure \ref*{fig:9a}}, there is a 5\% increase in the required contact reduction for scenario $\mathcal{M}$ to be replaced, as well as about a 1.5 day decrease between the required detection-and-removal day for the symptomatic individuals of the 4th cohort and the grouping of the 4th and 5th cohort. On the other hand, \hyperref[fig:13b]{Figure \ref*{fig:13b}} is identical to \hyperref[fig:9b]{Figure \ref*{fig:9b}}, meaning that the 5th cohort's contribution to the dynamics of $\mathcal{R}_0$ is minimal. This is further proved in \hyperref[fig:13c]{Figure \ref*{fig:13c}} and \hyperref[fig:13d]{Figure \ref*{fig:13d}}, where we see that the 2nd cohort dominates the dynamics of the simulation. In particular, when the contact reduction of the 2nd cohort is 50\%, scenario $\mathcal{L}$ can be replaced, whereas when the infectious individuals of the 2nd cohort are removed from the community at around the 4th day, scenario $\mathcal{M}$ can be replaced. 

\subsubsection{Epidemiological overview of the results, \texorpdfstring{$S_1$}{S1} versus \texorpdfstring{$S_2$}{S2}}
\label{Epidemiological overview of the results} 

Throughout our simulations we let $\left(a,b\right)$ of each gradation to take values in the 2D interval $[0,1)^2$ and illustrate the results in contour plots, where in the $x$-axis and $y$-axis we have $a\cdot 100\%$ and $b\cdot P_{I\to R}=\frac{b}{\gamma_I}=b\cdot 14$ days, respectively.

\begin{figure}[!h]
     \centering
        \begin{subfigure}[b]{0.48\textwidth}
            \centering
            \includegraphics[width=1.0\textwidth]{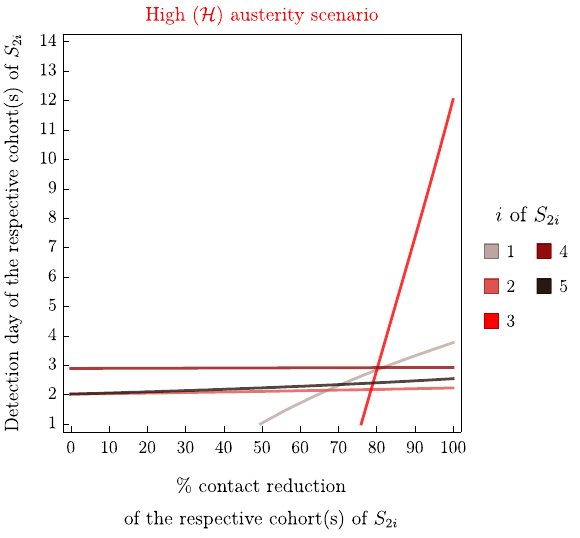}
            \caption{Loci of $\left(a\cdot 100\%, \frac{b}{\gamma_I} \right)$ of the respective substrategies of the family ${\left\{S_{2_i}\right\}}_{i=1}^{16}$ such that their $\mathcal{R}_0$ is equal to $0.571$, i.e, the $\mathcal{R}_0$ of scenario $\mathcal{H}$.}
            \label{fig:14a}
     \end{subfigure}
     \hfill
     \begin{subfigure}[b]{0.48\textwidth}
            \centering
            \includegraphics[width=1.0\textwidth]{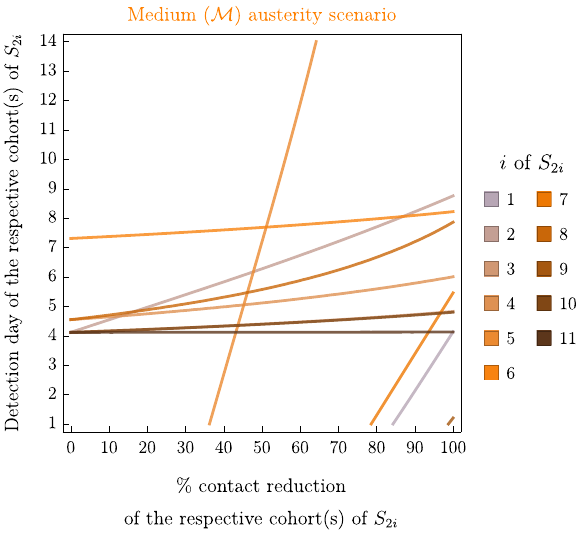}
            \caption{Loci of $\left(a\cdot 100\%, \frac{b}{\gamma_I} \right)$ of the respective substrategies of the family ${\left\{S_{2_i}\right\}}_{i=1}^{16}$ such that their $\mathcal{R}_0$ is equal to $1.427$, i.e, the $\mathcal{R}_0$ of scenario $\mathcal{M}$.}
            \label{fig:14b}
     \end{subfigure}
     \vskip\baselineskip
     \begin{subfigure}[b]{\textwidth}
            \centering
            \includegraphics[width=.48\textwidth]{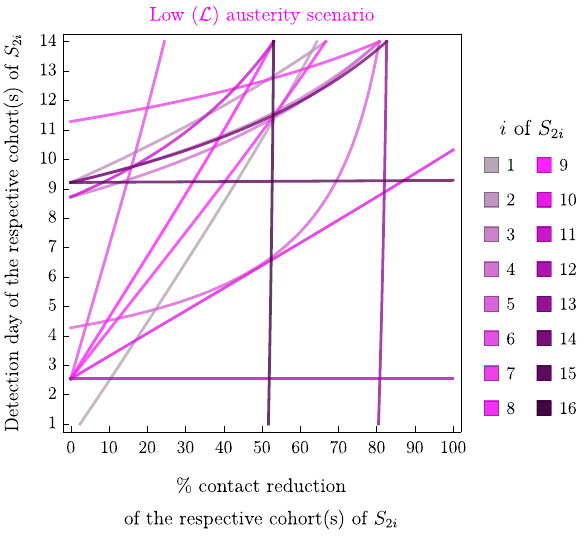}
            \caption{Loci of $\left(a\cdot 100\%, \frac{b}{\gamma_I} \right)$ of the respective substrategies of the family ${\left\{S_{2_i}\right\}}_{i=1}^{16}$ such that their $\mathcal{R}_0$ is equal to $2.283$, i.e, the $\mathcal{R}_0$ of scenario $\mathcal{L}$.}
            \label{fig:14c}
     \end{subfigure}

        \caption{Three contour plots illustrating the epidemiological coverage of each substrategy ($\mathcal{L, M}$ and $\mathcal{H}$) of horizontal lockdown strategy. The denser the plot is, the higher the epidemiological coverage.}
        \label{fig:14}
\end{figure}

\paragraph{$\mathcal{H}$ versus $S_2$} We begin by examining how many substrategies of ${\left\{S_{2_i}\right\}}_{i=1}^{16}$ can be considered as an alternative to scenario $\mathcal{H}$. As can be seen from \hyperref[fig:14a]{Figure \ref*{fig:14a}}, five substrategies of ${\left\{S_{2_i}\right\}}_{i=1}^{16}$ admit the same $\mathcal{R}_0$ as the respective one of scenario $\mathcal{H}$. Therefore, the epidemiological coverage of scenario $\mathcal{H}$ by the substrategies ${\left\{S_{2_i}\right\}}_{i=1}^{16}$ is 31.25\%. We highlight the fact that every one of the five substrategies that can replace scenario $\mathcal{M}$, regards restrictions on the 1st cohort.

\paragraph{$\mathcal{M}$ versus $S_2$} Next, we examine how many substrategies of ${\left\{S_{2_i}\right\}}_{i=1}^{16}$ can be considered as an alternative to scenario $\mathcal{M}$. As can be seen from \hyperref[fig:14b]{Figure \ref*{fig:14b}}, eleven substrategies of ${\left\{S_{2_i}\right\}}_{i=1}^{16}$ admit the same $\mathcal{R}_0$ as the respective one of scenario $\mathcal{M}$. Therefore, the epidemiological coverage of scenario $\mathcal{M}$ by the substrategies ${\left\{S_{2_i}\right\}}_{i=1}^{16}$ is 68.75\%.  We highlight the fact that every one of the eleven substrategies that can replace scenario $\mathcal{M}$, regards restrictions on the 1st and 2nd cohort.

\paragraph{$\mathcal{L}$ versus $S_2$} Finally, we examine how many substrategies of ${\left\{S_{2_i}\right\}}_{i=1}^{16}$ can be considered as an alternative to scenario $\mathcal{L}$. As can be seen from \hyperref[fig:14c]{Figure \ref*{fig:14c}}, all substrategies ${\left\{S_{2_i}\right\}}_{i=1}^{16}$ admit the same $\mathcal{R}_0$ as the respective one of scenario $\mathcal{L}$. Therefore, the epidemiological coverage of scenario $\mathcal{L}$ by the substrategies ${\left\{S_{2_i}\right\}}_{i=1}^{16}$ is 100\%.\\

A summary of the results of \hyperref[sec:VS]{\S \ref*{sec:VS}} can be seen in \hyperref[tab:AgeVSHorizontalComparison]{Table \ref*{tab:AgeVSHorizontalComparison}}.

\begin{table}[!h]
  \centering
  \resizebox{\linewidth}{!}{%
  \begin{tabular}{@{}ccccc@{}}
  \toprule
  \diagbox{Age-based restrictions}{Horizontal lockdowns}                                                         &  \textcolor{red}{High ($\mathcal{H}$)}                                                                            &    \textcolor{orange}{Medium ($\mathcal{M}$)}                                                                     &     \textcolor{magenta}{Low ($\mathcal{L}$)}                                                                       & \makecell{Epidemiological\\coverage}\\  \midrule
  \begin{tabular}[c]{@{}c@{}}Contact reduction: 1st, 2nd, 3rd cohorts \\  Testing: 4th, 5th cohorts\end{tabular} &  \begin{tabular}[c]{@{}c@{}}80.6\% and 3.06 days \\  88.\% and 6.4 days \\  95.4\% and 9.82 days\end{tabular}     &    \begin{tabular}[c]{@{}c@{}}42.3\% and 3.71 days \\  50.5\% and 7.43 days \\  58.9\% and 11.4 days\end{tabular} &     \begin{tabular}[c]{@{}c@{}}4.55\% and 11.7 days \\  12.5\% and 8.19 days \\  19.9\% and 4.59 days\end{tabular} & 100\%\\                                 \cmidrule(rl){1-1} \cmidrule(rl){2-4} \cmidrule(rl){5-5}
  \begin{tabular}[c]{@{}c@{}}Contact reduction: 4th, 5th cohorts \\  Testing: 1st, 2nd, 3rd cohorts\end{tabular} &  \begin{tabular}[c]{@{}c@{}}21.1\% and 2.92 days \\  50.4\% and 2.9 days \\  81.7\% and 2.89 days\end{tabular}    &    \begin{tabular}[c]{@{}c@{}}18.4\% and 7.43 days \\  50.\% and 7.68 days \\  78.9\% and 7.96 days\end{tabular}  &     \begin{tabular}[c]{@{}c@{}}17.8\% and 13.4 days \\  42.9\% and 12.4 days \\  67.9\% and 11.7 days\end{tabular} & 100\% \\                                \cmidrule(rl){1-1} \cmidrule(rl){2-4} \cmidrule(rl){5-5}
  \begin{tabular}[c]{@{}c@{}}Contact reduction: 1st cohort \\  Testing: 4th, 5th cohorts\end{tabular}            &  \xmark                                                                                                           &    \begin{tabular}[c]{@{}c@{}}98.8\% and 1 day \\  99.1\% and 1.05 days \\  99.5\% and 1.12 days\end{tabular}     &     \begin{tabular}[c]{@{}c@{}}10.7\% and 11.9 days \\  34.\% and 8.18 days \\  55.4\% and 4.31 days\end{tabular}  & 66.66\%\\                               \cmidrule(rl){1-1} \cmidrule(rl){2-4} \cmidrule(rl){5-5}
  \begin{tabular}[c]{@{}c@{}}Contact reduction: 4th, 5th cohorts \\  Testing: 1st cohort\end{tabular}            &  \begin{tabular}[c]{@{}c@{}}18.8\% and 2.05 days \\  48.6\% and 2.1 days \\  80.2\% and 2.18 days\end{tabular}    &    \begin{tabular}[c]{@{}c@{}}17.9\% and 4.72 days \\  50.\% and 5.11 days \\  79.1\% and 5.57 days\end{tabular}  &     \begin{tabular}[c]{@{}c@{}}17.5\% and 9.41 days \\  43.2\% and 10.8 days \\  68.6\% and 12.7 days\end{tabular} & 100\%\\                                 \cmidrule(rl){1-1} \cmidrule(rl){2-4} \cmidrule(rl){5-5}
  \begin{tabular}[c]{@{}c@{}}Contact reduction: 2nd cohort \\  Testing: 4th, 5th cohorts\end{tabular}            &  \xmark                                                                                                           &    \begin{tabular}[c]{@{}c@{}}82.7\% and 1.86 days \\  89.6\% and 3.29 days \\  96.4\% and 4.71 days\end{tabular} &     \begin{tabular}[c]{@{}c@{}}10.6\% and 11.7 days \\  27.7\% and 8.32 days \\  43.\% and 4.73 days\end{tabular}  & 66.66\%\\                                 \cmidrule(rl){1-1} \cmidrule(rl){2-4} \cmidrule(rl){5-5}
  \begin{tabular}[c]{@{}c@{}}Contact reduction: 4th, 5th cohorts \\  Testing: 2nd cohort\end{tabular}            &  \xmark                                                                                                           &    \begin{tabular}[c]{@{}c@{}}21.2\% and 4.65 days \\  50\% and 4.4 days \\  82.1\% and 4.22 days\end{tabular}    &     \begin{tabular}[c]{@{}c@{}}17.9\% and 9.82 days \\  43.1\% and 11 days \\  67.9\% and 12.7 days\end{tabular}   & 66.66\%\\                               \cmidrule(rl){1-1} \cmidrule(rl){2-4} \cmidrule(rl){5-5}
  \begin{tabular}[c]{@{}c@{}}Contact reduction: 3rd cohort \\  Testing: 4th, 5th cohorts\end{tabular}            &  \xmark                                                                                                           &    \xmark                                                                                                         &     \begin{tabular}[c]{@{}c@{}}19.5\% and 8.89 days \\  50.9\% and 6.45 days \\  82.1\% and 4.04 days\end{tabular} & 33.33\%\\                               \cmidrule(rl){1-1} \cmidrule(rl){2-4} \cmidrule(rl){5-5}
  \begin{tabular}[c]{@{}c@{}}Contact reduction: 4th, 5th cohorts \\  Testing: 3rd cohort\end{tabular}            &  \xmark                                                                                                           &    \xmark                                                                                                         &     \begin{tabular}[c]{@{}c@{}}23.1\% and 4.93 days \\  54.7\% and 6.89 days \\  75.\% and 11 days\end{tabular}    & 33.33\%\\                               \cmidrule(rl){1-1} \cmidrule(rl){2-4} \cmidrule(rl){5-5}
  \begin{tabular}[c]{@{}c@{}}Contact reduction: 1st cohort \\  Testing: 2nd cohort\end{tabular}                  &  \begin{tabular}[c]{@{}c@{}}58.1\% and 3.29 days \\  74.\% and 2.52 days \\  89.5\% and 1.59 days\end{tabular}    &    \begin{tabular}[c]{@{}c@{}}19.6\% and 4.95 days \\  50\% and 6.27 days \\  80.2\% and 7.71 days\end{tabular}   &     \begin{tabular}[c]{@{}c@{}}14\% and 10.1 days \\  34.3\% and 11.4 days \\  54.5\% and 13  days\end{tabular}    & 100\%\\                                 \cmidrule(rl){1-1} \cmidrule(rl){2-4} \cmidrule(rl){5-5}
  \begin{tabular}[c]{@{}c@{}}Contact reduction: 2nd cohort \\  Testing: 1st cohort\end{tabular}                  &  \begin{tabular}[c]{@{}c@{}}20.2\% and 2.41 days \\  51.2\% and 2.23 days \\  82.\% and 2.09 days\end{tabular}    &    \begin{tabular}[c]{@{}c@{}}20.1\% and 4.88 days \\  52.5\% and 5.64 days \\  82.1\% and 6.77 days\end{tabular} &     \begin{tabular}[c]{@{}c@{}}10.7\% and 9.36 days \\  28.2\% and 10.8 days \\  44.4\% and 12.6 days\end{tabular} & 100\% \\                                \cmidrule(rl){1-1} \cmidrule(rl){2-4} \cmidrule(rl){5-5}
  \begin{tabular}[c]{@{}c@{}}Contact reduction: 4th cohort \\  Testing: 5th cohort\end{tabular}                  &  \xmark                                                                                                           &    \xmark                                                                                                         &     \begin{tabular}[c]{@{}c@{}}81.1\% and 3.46 days \\  81.8\% and 7.5 days \\  82.3\% and 11.2 days\end{tabular}  & 33.33\%\\                               \cmidrule(rl){1-1} \cmidrule(rl){2-4} \cmidrule(rl){5-5}
  \begin{tabular}[c]{@{}c@{}}Contact reduction: 5th cohort \\  Testing: 4th cohort\end{tabular}                  &  \xmark                                                                                                           &    \xmark                                                                                                         &     \begin{tabular}[c]{@{}c@{}}18.9\% and 2.53 days \\  48.9\% and 2.53 days \\  79.7\% and 2.53 days\end{tabular} & 33.33\%\\                               \cmidrule(rl){1-1} \cmidrule(rl){2-4} \cmidrule(rl){5-5}
  \begin{tabular}[c]{@{}c@{}}Contact reduction: 2nd cohort \\  Testing: 4th cohort\end{tabular}                  &  \xmark                                                                                                           &    \begin{tabular}[c]{@{}c@{}}87.5\% and 1.62 days \\  92.\% and 2.51 days \\  96.4\% and 3.41 days\end{tabular}  &     \begin{tabular}[c]{@{}c@{}}11.5\% and 2.76 days \\  33.8\% and 7.25 days \\  52.8\% and 11.3 days\end{tabular} & 66.66\%\\                               \cmidrule(rl){1-1} \cmidrule(rl){2-4} \cmidrule(rl){5-5}
  \begin{tabular}[c]{@{}c@{}}Contact reduction: 4th cohort \\  Testing: 2nd cohort\end{tabular}                  &  \xmark                                                                                                           &    \begin{tabular}[c]{@{}c@{}}21.4\% and 4.63 days \\  51.1\% and 4.4 days \\  81.5\% and 4.22 days\end{tabular}  &     \begin{tabular}[c]{@{}c@{}}17.9\% and 9.81 days \\  44.5\% and 11  days \\  68.1\% and 12.6 days\end{tabular}  & 66.66\%\\                               \cmidrule(rl){1-1} \cmidrule(rl){2-4} \cmidrule(rl){5-5}
  \begin{tabular}[c]{@{}c@{}}Contact reduction: 2nd cohort \\  Testing: 5th cohort\end{tabular}                  &  \xmark                                                                                                           &    \xmark                                                                                                         &     \begin{tabular}[c]{@{}c@{}}52.1\% and 3.52 days \\  52.5\% and 7.36 days \\  52.9\% and 11.3 days\end{tabular} & 33.33\%\\                               \cmidrule(rl){1-1} \cmidrule(rl){2-4} \cmidrule(rl){5-5}
  \begin{tabular}[c]{@{}c@{}}Contact reduction: 5th cohort \\  Testing: 2nd cohort\end{tabular}                  &  \xmark                                                                                                           &    \begin{tabular}[c]{@{}c@{}}20.4\% and 4.12 days \\  51.\% and 4.12 days \\  81.2\% and 4.12 days\end{tabular}  &     \begin{tabular}[c]{@{}c@{}}18.9\% and 9.22 days \\  49.1\% and 9.25 days \\  79.3\% and 9.27 days\end{tabular} & 66.66\%\\                               \cmidrule(rl){1-1} \cmidrule(rl){2-4} \cmidrule(rl){5-5}
    Social coverage                                                                                              &  31.25\%                                                                                                          & 68.75\%                                                                                                           &      100\%                                                                                                         & 66.66\% \\                              \bottomrule
  \end{tabular}
  }
  \caption{Horizontal lockdowns versus age-based restrictions. The total coverage of horizontal lockdowns from age-based restrictions is 66.66\%. Additionally, the table is populated with representative values of $\left(a\cdot 100\%, \frac{b}{\gamma_I} \right)$ of the strategic scale that each age-based strategy needs to have in order for the strategy to have the same $\mathcal{R}_0$ as each of the three horizontal lockdown scenarios.}
  \label{tab:AgeVSHorizontalComparison}
  \end{table}




\clearpage

\afterpage{\clearpage} 

\section{Conclusions and discussion}
\label{sec:CD}
In this paper, we introduced a scheme for the comparison of certain types of interventions for the restriction of an epidemiological phenomenon. This scheme incorporates some novel notions such as \textquote{strategy}, \textquote{substrategy}, \textquote{gradable strategy} and its \textquote{gradation}, \textquote{comparison table}, as well as \textquote{epidemiological coverage} and \textquote{social coverage}. Then, we utilized the aforementioned scheme and the age-based epidemiological compartment problem  studied in \cite{bgt2023SVEAIRmodel} to compare horizontal lockdown policies with various age-based interventions. 

In particular, we distributed the total population into five cohorts, based on the age of each individual (in ascending order) and we defined the graded strategy of horizontal lockdowns, considering three scenarios of horizontal lockdowns with varying intensity, Low ($\mathcal{L}$), Medium ($\mathcal{M}$) and High ($\mathcal{H}$). We also defined the strategy of age-based restrictions, consisting of $16$ substrategies. In general, our results suggest that these two strategies are comparable mainly at low or medium level of intensity. Particularly, throughout our simulations, which used data from the literature, we deduced that the strategies that targeted the 1st and 2nd cohort had the best epidemiological coverage. Moreover, all substrategies were able to admit the same $\mathcal{R}_0$ as the respective one of scenario $\mathcal{L}$, meaning a 100\% social coverage of $\mathcal{L}$, while the social coverage of scenarios $\mathcal{M}$ and $\mathcal{H}$ by the substrategies is 68.75\% and 31.25\%, respectively. 

Future work could entail the generalization of the notion of strategy, hence the comparison process itself.

\section*{Declaration of interests}
The authors declare that they have no known competing financial interests or personal relationships that could have appeared to influence the work reported in this paper.

\begin{appendices} \label{appendix}

\section{The employed epidemiological model}
\label{sec:model} 

Here we use the epidemiological model $\mathscr{M}$ along with the respective problem $\mathscr{P}$, introduced and studied in \cite{bgt2023SVEAIRmodel}, as a means of utilization of the proposed scheme in answering the main question of the present paper. 
We choose this model as it incorporates both symptomatic and asymptomatic infectious individuals, with the latter playing an important role in the spread of the COVID-19 (see \cite{gao2021systematic} and many references therein), as well as the age of the infected/infectious individuals. 

After \textit{scaling} the independent age-variable, $\theta$, and turning it to another time-variable measured in the same units as $t$ (see \cite{bgt2023SVEAIRmodel}) and using the relation $N=S+V+E+A+I+R$, we obtain the following model


\begin{subequations}
\label{SVEIAR-age-scl}
\begin{align}
&\begin{cases}
\dfrac{\de S}{\de t}=\mu N_0-\left(p+\int\limits_0^\infty{\beta_A{\left(\theta\right)}a{\left(\,\cdot\,,\theta\right)}+\beta_I{\left(\theta\right)}i{\left(\,\cdot\,,\theta\right)}\,\de\theta}+\mu\right)S\\
S{\left(0\right)}=S_0,\label{SVEIAR-age-scl;a}
\end{cases}\\
&\begin{cases}
\dfrac{\de V}{\de t}=pS-\left(\zeta\epsilon+\int\limits_0^\infty{\beta_A{\left(\theta\right)}a{\left(\,\cdot\,,\theta\right)}+\beta_I{\left(\theta\right)}i{\left(\,\cdot\,,\theta\right)}\,\de\theta}\left(1-\epsilon\right)+\mu\right)V\\
V{\left(0\right)}=V_0,\label{SVEIAR-age-scl;b}
\end{cases}\\
&\begin{cases}
\dfrac{\partial e}{\partial t}+\dfrac{\partial e}{\partial \theta}=-\left(k+\mu\right)e\\
e{\left(\,\cdot\,,0\right)}=\int\limits_0^\infty{\beta_A{\left(\theta\right)}a{\left(\,\cdot\,,\theta\right)}+\beta_I{\left(\theta\right)}i{\left(\,\cdot\,,\theta\right)}\,\de\theta}\left(S+\left(1-\epsilon\right)V\right)\\
e{\left(0,\,\cdot\,\right)}=e_0,\label{SVEIAR-age-scl;c}
\end{cases}\\
&\begin{cases}
\dfrac{\partial a}{\partial t}+\dfrac{\partial a}{\partial \theta}=-\left(\gamma_A\xi+\chi\left(1-\xi\right)+\mu\right)a\\
a{\left(\,\cdot\,,0\right)}=\int\limits_0^\infty{k{\left(\theta\right)}q{\left(\theta\right)}e{\left(\,\cdot\,,\theta\right)}\,\de\theta}\\
a{\left(0,\,\cdot\,\right)}=a_0,\label{SVEIAR-age-scl;d}
\end{cases}\\
&\begin{cases}
\dfrac{\partial i}{\partial t}+\dfrac{\partial i}{\partial \theta}=-\left(\gamma_I+\mu\right)i\\
i{\left(\,\cdot\,,0\right)}=\int\limits_0^\infty{k{\left(\theta\right)}\left(1-q{\left(\theta\right)}\right)e{\left(\,\cdot\,,\theta\right)}+\chi{\left(\theta\right)}\left(1-\xi{\left(\theta\right)}\right)a{\left(\,\cdot\,,\theta\right)}\,\de\theta}\\
i{\left(0,\,\cdot\,\right)}=i_0.\label{SVEIAR-age-scl;e}
\end{cases}
\end{align}
\end{subequations}

\begin{figure}[!h]
\centering
\includegraphics[width=0.9\textwidth]{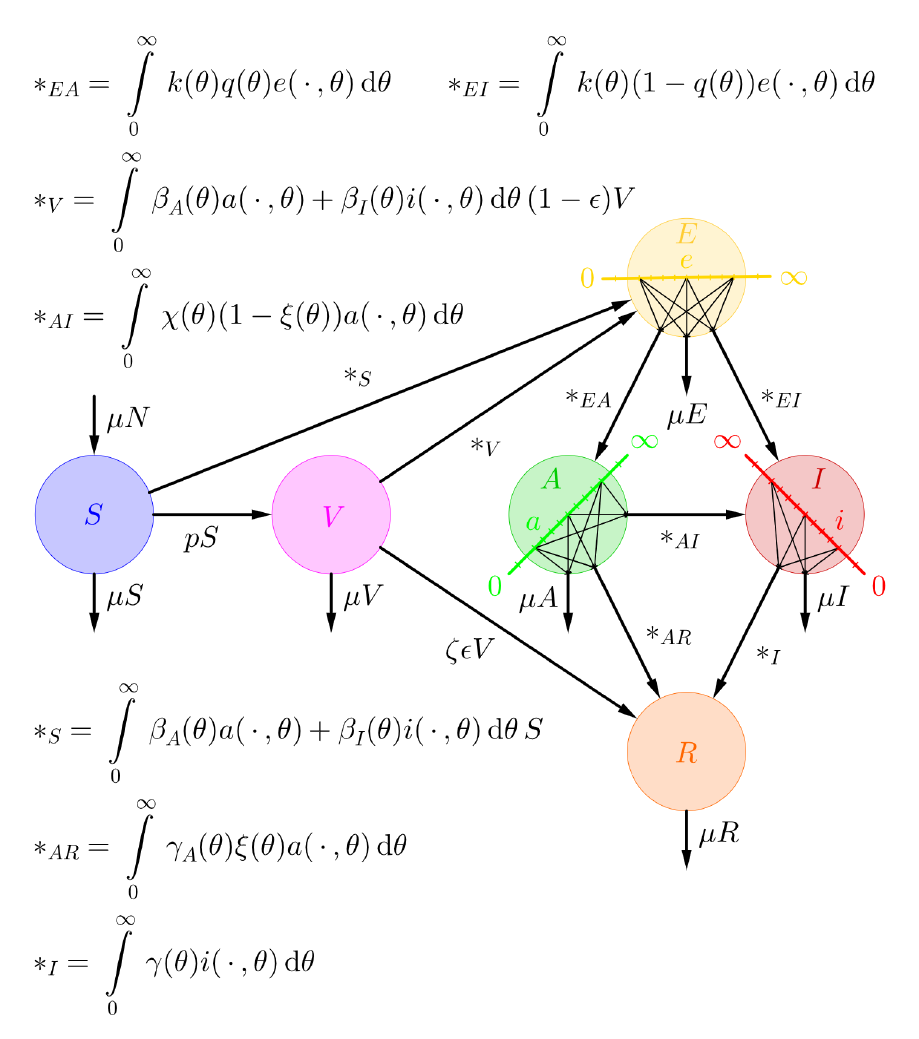}
\caption{Flow diagram of $\mathscr{P}$ \eqref{SVEIAR-age-scl}.}
\label{fig:2}
\end{figure}

The flow diagram of the differential equations in \eqref{SVEIAR-age-scl} is shown in \hyperref[fig:2]{Figure \ref*{fig:2}}, and the dimensional units of all variables and parameters appeared in $\mathscr{P}$ \eqref{SVEIAR-age-scl} are gathered in \hyperref[table-1]{Table \ref*{table-1}}.

\begin{table}[p]
\centering 
\begin{tabular}{p{2cm}p{11cm}p{2cm}}
\hline
Independent variables  & Description  & Units  \\ [0.5ex]
\hline         
$t$ & Time   & T \\
$\theta$ & Age, i.e., time elapsed since, e.g., birth or infection & $\Theta$ \\  [1ex]      
\hline
\hline
Conversion factor  & Description  & Units  \\ [0.5ex]
\hline
$\omega$ & Conversion factor from the units of $\theta$ to the units of $t$ & T$\,\Theta^{-1}$\\  [1ex] 
\hline
\hline
Dependent variables  & Description  & Units  \\ [0.5ex]
\hline        
$N$ & Number of total population of individuals  & \#  \\
$S$ & Number of susceptible individuals   & \#  \\
$V$ & Number of vaccinated-with-a-prophylactic-vaccine individuals   & \#  \\
$e$ & Age density of latent/exposed individuals  &  \#$\,\Theta^{-1}$ \\
$E$ & Number of latent/exposed individuals   & \#  \\
$a$ & Age density of asymptomatic infectious individuals   &  \#$\,\Theta^{-1}$  \\
$A$ & Number of asymptomatic infectious individuals   & \#  \\
$i$ & Age density of symptomatic infectious   &   \#$\,\Theta^{-1}$ \\
$I$ & Number of symptomatic infectious individuals   & \#  \\
$R$ & Number of recovered/removed individuals   & \#  \\  [1ex]      
\hline
\hline
Parameters & Description  & Units  \\ [0.5ex]
\hline
$N_0$ & Population size & \# \\
$\mu$ & Birth/Death rate   & T$^{-1}$ \\
$\beta_A$ & Transmission rate of asymptomatic infectious individuals & \#$^{-1}\,$T$^{-1}$ \\
$\beta_I$ & Transmission rate of symptomatic infectious individuals & \#$^{-1}\,$T$^{-1}$ \\
$p$ & Vaccination rate &  T$^{-1}$ \\
$\epsilon$ & Vaccine effectiveness  & -  \\
$\zeta$ & Vaccine-induced immunity rate & T$^{-1}$  \\
$k$ & Latent rate (rate of susceptible individuals becoming infectious)  & T$^{-1}$ \\
$q$ & Proportion of the latent/exposed individuals becoming asymptomatic infectious  & - \\
$\xi$ &  Proportion of the asymptomatic infectious individuals becoming recovered/removed (without developing any symptoms) & - \\
$\chi$ &  Incubation rate (rate of a part of asymptomatic infectious individuals developing symptoms) & T$^{-1}$ \\
$\gamma_A$ & Recovery rate of asymptomatic infectious individuals  &  T$^{-1}$\\
$\gamma_I$ & Recovery rate of symptomatic infectious individuals  &  T$^{-1}$ \\ [1ex]      
\hline
\end{tabular}
\caption{Description of the independent and dependent variables and parameters of $\mathscr{M}$, along with their units.}
\label{table-1}
\end{table}


From the analysis conducted in \cite{bgt2023SVEAIRmodel}, the \textit{basic reproductive number}, $\mathcal{R}_0$, of the model is
\begin{equation}
    \boxed{\mathbb{R}_0^+\ni\mathcal{R}_0\coloneqq\frac{\mu N_0}{p+\mu}\left(1+\frac{p\left(1-\epsilon\right)}{\zeta\epsilon+\mu}\right)\left(\mathcal{R}_A+\mathcal{R}_I\right),}
    \label{R0-definition}
\end{equation}
where 
\begin{equation*}
\boxed{\mathbb{R}_0^+\ni\mathcal{R}_A\coloneqq\int\limits_0^\infty{k{\left(s\right)}q{\left(s\right)}\e^{-\int\limits_0^s{k{\left(\tau\right)}+\mu\,\de\tau}}\,\de s}\int\limits_0^\infty{\beta_A{\left(s\right)}\e^{-\int\limits_0^s{\gamma_A{\left(\tau\right)}\xi{\left(\tau\right)}+\chi{\left(\tau\right)}\left(1-\xi{\left(\tau\right)}\right)+\mu\,\de\tau}}\,\de s}}
\end{equation*}
and 
\begin{equation*}
\boxed{\begin{aligned}
\mathbb{R}_0^+\ni\mathcal{R}_I&\coloneqq\left(
\begin{aligned}
&\phantom{+}\int\limits_0^\infty{k{\left(s\right)}\left(1-q{\left(s\right)}\right)\e^{-\int\limits_0^s{k{\left(\tau\right)}+\mu\,\de\tau}}\,\de s}+\\
&+\int\limits_0^\infty{k{\left(s\right)}q{\left(s\right)}\e^{-\int\limits_0^s{k{\left(\tau\right)+\mu}\,\de\tau}}\,\de s}\int\limits_0^\infty{\chi{\left(s\right)}\left(1-\xi{\left(s\right)}\right)\e^{-\int\limits_0^s{\gamma_A{\left(\tau\right)}\xi{\left(\tau\right)}+\chi{\left(\tau\right)}\left(1-\xi{\left(\tau\right)}\right)+\mu\,\de\tau}}\,\de s}
\end{aligned}
\right)\times\\
&\phantom{\coloneqq}\times\int\limits_0^\infty{\beta_I{\left(s\right)}\e^{-\int\limits_0^s{\gamma_I{\left(\tau\right)}+\mu\,\de\tau}}\,\de s}.
\end{aligned}}
\end{equation*}

\section{Parameter estimation} \label{sec:paramestim}

We now present parameter values fitting for the case of SARS-CoV-2. The chosen values are taken from the biological and medical literature. Bellow, we give a detailed explanation about the value of each parameter, whereas a summary of the parameter values can be found in \hyperref[tab:paramValues]{Table \ref*{tab:paramValues}}.

The size of the population, $N_0 = 80 \cdot 10^6$ individuals, is assumed to be that of a relative large country, such as Germany, Turkey, or Thailand \citep{ourworldindata}.

The birth/death rate, $\mu = 4.38356 \cdot 10^{-5} \; \text{day}^{-1}$, is taken from data from \cite{ourworldindata}. The average birth/death rate of the world for the year 2021 is about 16 per 1000 individuals per year. Hence, we convert the aforementioned quantity from ``per 1000 individuals per year" to ``per day" to get
\begin{equation*}
    16 \frac{1}{1000 \; \text{individuals} \cdot \text{year}} \mapsto 16 \cdot 10^{-3} \frac{1}{365 \; \text{days}} = 4.38356 \cdot 10^{-5} \; \text{day}^{-1}=\mu.
\end{equation*}

\begin{figure}[h]
    \includegraphics[width=\textwidth]{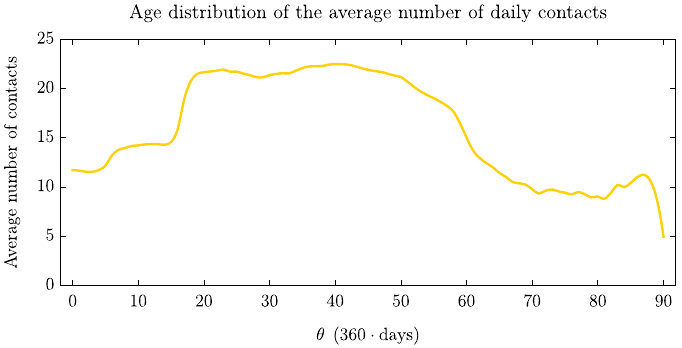}
    \caption{Age density (in years) of the average number of daily contacts, $c$, taken from \cite{DelValle_Hyman_Hethcote_Eubank_2007}.}
    \label{fig:3}
\end{figure}

\begin{figure}[ht]
    \includegraphics[width=\textwidth]{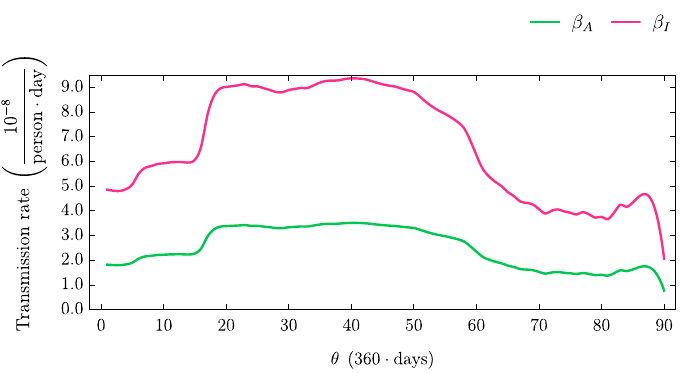}
    \caption{Estimation of the age density (in years) of the asymptomatic and symptomatic transmission rate, assuming $\beta_A = \frac{c_A \cdot \varpi_{E\to A}}{N_0}$ and $ \beta_I= \frac{c_I \cdot \varpi_{E\to I}}{N_0}$ according to \eqref{beta-def}, where $c_A=c=c_I$.}
    \label{fig:4}
\end{figure}

For the transmission rate of asymptomatic and symptomatic infectious individuals, we firstly assume the probability of an exposed individual passing to the compartments of asymptomatic and symptomatic individuals to be $\varpi_{E\to A} = \frac{1}{8}$ and $\varpi_{E\to I} = \frac{1}{3}$, respectively. From \cite{DelValle_Hyman_Hethcote_Eubank_2007}, we have that the average number of daily contacts of any person, regardless its epidemiological status, of age $\theta$, $c(\theta)$, follows the graph as seen in \hyperref[fig:3]{Figure \ref*{fig:3}}. To digitize the data of the contacts, we use WebPlotDigitizer 4.6 \citep{Rohatgi2022} to manually extract data points from Fig. 2 of \cite{DelValle_Hyman_Hethcote_Eubank_2007} and then interpolated them using a third order polynomial interpolation scheme through Mathematica 13.1 \citep{Mathematica} and the function \texttt{Interpolation}. Subsequently, from \eqref{beta-def} we deduce that the transmission rate of asymptomatic and symptomatic infectious individuals are the functions presented in \hyperref[fig:4]{Figure \ref*{fig:4}}. 

The vaccination rate, $p = 10^{-3} \; \text{day}^{-1}$, is taken from data from \cite{ourworldindata}, during the summer of 2021 in the USA, when the Delta variant of SARS-CoV-2 was the dominant variant. During the end of summer, the percentage of fully vaccinated USA citizens was about 54\% whereas in the beginning of summer it was around 45\%. Hence, we estimate the vaccination from that three-month period to be $p = \frac{54\% - 45\%}{90 }\; \text{day}^{-1} = 10^{-3} \; \text{day}^{-1}.$

The vaccine effectiveness, $\epsilon = 0.7$, is estimated from data from \cite{grant2022impact}. In \cite{grant2022impact}, the authors find that with the BNT162b2 vaccine, the effectiveness of two doses is 88.0\% among those with the Delta variant, whereas with the ChAdOx1 nCoV-19 vaccine, the respective effectiveness of two doses was 67.0\%. Hence, we assume $\epsilon = 0.7$.

The vaccine-induced immunity rate, $\zeta = \frac{1}{14} \; \text{day}^{-1}$, is estimated from data from \cite{chau2022immunogenicity}. The authors of \cite{chau2022immunogenicity} report that, after two weeks of the second dose of the ChAdOx1 nCoV-19 vaccine, the percentage of study participants with detectable neutralizing antibodies reached 98.1\%. 

The latent rate, $k$, is found by estimating that the latent and incubation period differ by 1 day. In \cite{kang2022transmission}, the authors examined data from 93 Delta transmission pairs and estimated the latent period by fitting the data to the Weibull distribution, which made the best fit. They found the mean latent period to be 3.9 days. In \cite{wu2022incubation}, the authors performed a systematic review and meta-analysis of 141 articles and found that the incubation periods of COVID-19 caused by the Alpha, Beta, Delta, and Omicron variants were 5.00, 4.50, 4.41, and 3.42 days, respectively. Hence, assuming that the latent and incubation period vary by 1 day, we have that $k = \frac{\chi}{1-\chi}$, and by substituting $\chi$ as found later in the present section, we have that
\begin{equation} \label{eq:valuek}
    k(\theta) = \begin{cases}
    \frac{1}{4}   \; \text{day}^{-1},   & \theta < 30 \cdot 360 \; \text{day}\\
    \frac{1}{4.8} \; \text{day}^{-1},   & 30 \cdot 360 \; \text{day} \; \le \theta < 40 \cdot 360 \; \text{day}\\
    \frac{1}{4.8} \; \text{day}^{-1},   & 40 \cdot 360 \; \text{day} \; \le \theta < 50 \cdot 360 \; \text{day}\\
    \frac{1}{5.5} \; \text{day}^{-1},   & 50 \cdot 360 \; \text{day} \; \le \theta < 60 \cdot 360 \; \text{day}\\
    \frac{1}{3.1} \; \text{day}^{-1},   & 60 \cdot 360 \; \text{day} \; \le \theta < 70 \cdot 360 \; \text{day}\\
    \frac{1}{6}   \; \text{day}^{-1},   & 70 \cdot 360 \; \text{day} \; \le \theta \; .
    \end{cases}
\end{equation}
where $\theta$ is measured in years.

\begin{figure}[!h]
    \centering
    \includegraphics[width=\textwidth]{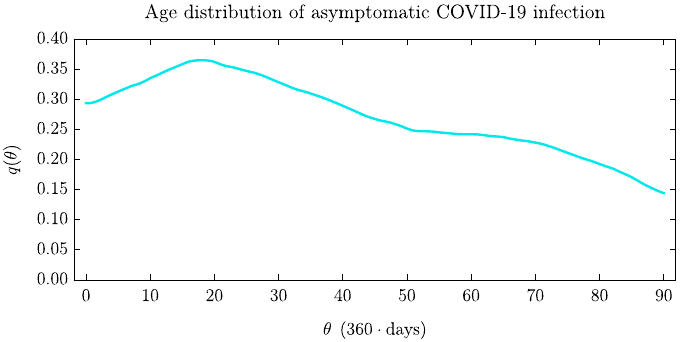}
    \caption{Percentage of asymptomatic COVID-19 infection, by age (in years), taken from \cite{sah2021asymptomatic}.}
    \label{fig:5}
\end{figure}

The proportion of the latent/exposed individuals becoming asymptomatic
infectious, $q$, is taken from \cite{sah2021asymptomatic}, where the authors estimated the asymptomatic proportion by age, by performing a systematic review and meta-analysis of 38 studies involving 14850 individuals. The curve they estimated can be seen in \hyperref[fig:5]{Figure \ref*{fig:5}}. To digitize the data, we use the same procedure we used for the age density of daily contacts described earlier in the present section.

The proportion of the asymptomatic infectious individuals becoming recovered/removed without developing any symptoms, $\xi = 0.5$, is estimated from data from \cite{he2021proportion} and \cite{buitrago2022occurrence}. In \cite{he2021proportion}, the authors performed a systematic review and meta-analysis from 41 studies containing the rate of asymptomatic COVID-19 infection before May 20, 2020, aggregating 50155 patients, and found that nearly half of the patients with no symptoms at the time of their detection, would develop symptoms later. In \cite{buitrago2022occurrence}, the authors performed a systematic review and meta-analysis from 130 studies and reported the percentage of persistently asymptomatic individuals being between 14 to 50\%. Hence, we choose the proportion of persistently asymptomatic individuals being 50\%.

\afterpage{\clearpage}

The incubation rate, $\chi$, is estimated from data from \cite{tan2020does}. The authors of \cite{tan2020does} found that the incubation period varies with age, based on data from Singaporean hospitals between January 23, 2020 and April 2, 2020. The authors divided the participants based on their age (in years) to six groups ($<$30, 30--39, 40--49, 50--59, 60--69 and 70$<$) and presented their results through a box plot. Hence, we assume that $\chi$ is a piecewise function with its domain intervals being the six aforementioned age groups, and with the function being constant on each interval and equal to one over the median of the respective age group. Therefore, we have that
\begin{equation}  \label{eq:valuechi}
    \chi(\theta) = \begin{cases}
    \frac{1}{5}   \; \text{day}^{-1},   & \theta < 30 \cdot 360 \; \text{day}\\
    \frac{1}{5.8} \; \text{day}^{-1},   & 30 \cdot 360 \; \text{day} \; \le \theta < 40 \cdot 360 \; \text{day}\\
    \frac{1}{5.8} \; \text{day}^{-1},   & 40 \cdot 360 \; \text{day} \; \le \theta < 50 \cdot 360 \; \text{day}\\
    \frac{1}{6.5} \; \text{day}^{-1},   & 50 \cdot 360 \; \text{day} \; \le \theta < 60 \cdot 360 \; \text{day}\\
    \frac{1}{4.1} \; \text{day}^{-1},   & 60 \cdot 360 \; \text{day} \; \le \theta < 70 \cdot 360 \; \text{day}\\
    \frac{1}{7}   \; \text{day}^{-1},   & 70 \cdot 360 \; \text{day} \; \le \theta \; .
    \end{cases}
\end{equation}
where $\theta$ is measured in years.

The recovery rate of asymptomatic infectious individuals, $\gamma_A = \frac{1}{8} \; \text{day}^{-1}$ , and recovery rate of the symptomatic infectious individuals, $\gamma_I = \frac{1}{14} \; \text{day}^{-1}$, is estimated from \cite{byrne2020inferred}. In \cite{byrne2020inferred}, the authors performed a rapid scoping review up to April 1, 2020 and found that the median infectious period for asymptomatic cases was 6.5--9.5 days, whereas time from symptom onset to two negative RT-PCR tests ranged from 10.9 to 15.8 days. Hence, we assume that the recovery period of asymptomatic and symptomatic infectious individuals to be 8 and 14 days respectively.


\end{appendices}

\bibliographystyle{elsarticle-num-names} 
\bibliography{mybibfile} 

\begin{thebibliography}{55}
\expandafter\ifx\csname natexlab\endcsname\relax\def\natexlab#1{#1}\fi
\providecommand{\url}[1]{\texttt{#1}}
\providecommand{\href}[2]{#2}
\providecommand{\path}[1]{#1}
\providecommand{\DOIprefix}{doi:}
\providecommand{\ArXivprefix}{arXiv:}
\providecommand{\URLprefix}{URL: }
\providecommand{\Pubmedprefix}{pmid:}
\providecommand{\doi}[1]{\href{http://dx.doi.org/#1}{\path{#1}}}
\providecommand{\Pubmed}[1]{\href{pmid:#1}{\path{#1}}}
\providecommand{\bibinfo}[2]{#2}
\ifx\xfnm\relax \def\xfnm[#1]{\unskip,\space#1}\fi
\bibitem[{Bitsouni et~al.(2024)Bitsouni, Gialelis, and Tsilidis}]{bgt2023SVEAIRmodel}
\bibinfo{author}{V.~Bitsouni}, \bibinfo{author}{N.~Gialelis}, \bibinfo{author}{V.~Tsilidis},
\newblock \bibinfo{title}{An age-structured {SVEAIR} epidemiological model},
\newblock \bibinfo{journal}{Mathematical Methods in the Applied Sciences}  (\bibinfo{year}{2024}). \DOIprefix\doi{10.1002/mma.10165}.
\bibitem[{Brodeur et~al.(2021)Brodeur, Gray, Islam, and Bhuiyan}]{brodeur2021literature}
\bibinfo{author}{A.~Brodeur}, \bibinfo{author}{D.~Gray}, \bibinfo{author}{A.~Islam}, \bibinfo{author}{S.~Bhuiyan},
\newblock \bibinfo{title}{A literature review of the economics of {COVID-19}},
\newblock \bibinfo{journal}{Journal of Economic Surveys} \bibinfo{volume}{35} (\bibinfo{year}{2021}) \bibinfo{pages}{1007--1044}. \DOIprefix\doi{10.1111/joes.12423}.
\bibitem[{Chen et~al.(2021)Chen, Vullikanti, Santos, Venkatramanan, Hoops, Mortveit, Lewis, You, Eubank, Marathe et~al.}]{chen2021epidemiological}
\bibinfo{author}{J.~Chen}, \bibinfo{author}{A.~Vullikanti}, \bibinfo{author}{J.~Santos}, \bibinfo{author}{S.~Venkatramanan}, \bibinfo{author}{S.~Hoops}, \bibinfo{author}{H.~Mortveit}, \bibinfo{author}{B.~Lewis}, \bibinfo{author}{W.~You}, \bibinfo{author}{S.~Eubank}, \bibinfo{author}{M.~Marathe}, et~al.,
\newblock \bibinfo{title}{Epidemiological and economic impact of {COVID-19} in the {US}},
\newblock \bibinfo{journal}{Scientific Reports} \bibinfo{volume}{11} (\bibinfo{year}{2021}) \bibinfo{pages}{20451}. \DOIprefix\doi{10.1038/s41598-021-99712-z}.
\bibitem[{Deb et~al.(2021)Deb, Furceri, Ostry, and Tawk}]{Deb2021}
\bibinfo{author}{P.~Deb}, \bibinfo{author}{D.~Furceri}, \bibinfo{author}{J.~D. Ostry}, \bibinfo{author}{N.~Tawk},
\newblock \bibinfo{title}{The economic effects of {COVID}-19 containment measures},
\newblock \bibinfo{journal}{Open Economies Review} \bibinfo{volume}{33} (\bibinfo{year}{2021}) \bibinfo{pages}{1--32}. \DOIprefix\doi{10.1007/s11079-021-09638-2}.
\bibitem[{Mathieu et~al.(2020)Mathieu, Ritchie, Rodés-Guirao, Appel, Giattino, Hasell, Macdonald, Dattani, Beltekian, Ortiz-Ospina, and Roser}]{ourworldindata}
\bibinfo{author}{E.~Mathieu}, \bibinfo{author}{H.~Ritchie}, \bibinfo{author}{L.~Rodés-Guirao}, \bibinfo{author}{C.~Appel}, \bibinfo{author}{C.~Giattino}, \bibinfo{author}{J.~Hasell}, \bibinfo{author}{B.~Macdonald}, \bibinfo{author}{S.~Dattani}, \bibinfo{author}{D.~Beltekian}, \bibinfo{author}{E.~Ortiz-Ospina}, \bibinfo{author}{M.~Roser},
\newblock \bibinfo{title}{{Coronavirus Pandemic (COVID-19)}},
\newblock \bibinfo{journal}{Our World in Data}  (\bibinfo{year}{2020}). \bibinfo{note}{\url{https://ourworldindata.org/coronavirus}}.
\bibitem[{Asahi et~al.(2021)Asahi, Undurraga, Vald{\'e}s, and Wagner}]{asahi2021effect}
\bibinfo{author}{K.~Asahi}, \bibinfo{author}{E.~A. Undurraga}, \bibinfo{author}{R.~Vald{\'e}s}, \bibinfo{author}{R.~Wagner},
\newblock \bibinfo{title}{The effect of {COVID-19} on the economy: Evidence from an early adopter of localized lockdowns},
\newblock \bibinfo{journal}{Journal of Global Health} \bibinfo{volume}{11} (\bibinfo{year}{2021}). \DOIprefix\doi{10.7189/jogh.10.05002}.
\bibitem[{Karatayev et~al.(2020)Karatayev, Anand, and Bauch}]{karatayev2020local}
\bibinfo{author}{V.~A. Karatayev}, \bibinfo{author}{M.~Anand}, \bibinfo{author}{C.~T. Bauch},
\newblock \bibinfo{title}{Local lockdowns outperform global lockdown on the far side of the {COVID-19} epidemic curve},
\newblock \bibinfo{journal}{Proceedings of the National Academy of Sciences} \bibinfo{volume}{117} (\bibinfo{year}{2020}) \bibinfo{pages}{24575--24580}. \DOIprefix\doi{10.1073/pnas.2014385117}.
\bibitem[{Perra(2021)}]{perra2021non}
\bibinfo{author}{N.~Perra},
\newblock \bibinfo{title}{Non-pharmaceutical interventions during the {COVID-19} pandemic: A review},
\newblock \bibinfo{journal}{Physics Reports} \bibinfo{volume}{913} (\bibinfo{year}{2021}) \bibinfo{pages}{1--52}. \DOIprefix\doi{10.1016/j.physrep.2021.02.001}.
\bibitem[{Demers et~al.(2023)Demers, Fagan, Potluri, and Calabrese}]{demers2023relationship}
\bibinfo{author}{J.~Demers}, \bibinfo{author}{W.~F. Fagan}, \bibinfo{author}{S.~Potluri}, \bibinfo{author}{J.~M. Calabrese},
\newblock \bibinfo{title}{{The relationship between controllability, optimal testing resource allocation, and incubation-latent period mismatch as revealed by COVID-19}},
\newblock \bibinfo{journal}{Infectious Disease Modelling} \bibinfo{volume}{8} (\bibinfo{year}{2023}) \bibinfo{pages}{514--538}. \DOIprefix\doi{10.1016/j.idm.2023.04.007}.
\bibitem[{Adegbite et~al.(2023)Adegbite, Edeki, Isewon, Emmanuel, Dokunmu, Rotimi, Oyelade, and Adebiyi}]{adegbite2023mathematical}
\bibinfo{author}{G.~Adegbite}, \bibinfo{author}{S.~Edeki}, \bibinfo{author}{I.~Isewon}, \bibinfo{author}{J.~Emmanuel}, \bibinfo{author}{T.~Dokunmu}, \bibinfo{author}{S.~Rotimi}, \bibinfo{author}{J.~Oyelade}, \bibinfo{author}{E.~Adebiyi},
\newblock \bibinfo{title}{Mathematical modeling of malaria transmission dynamics in humans with mobility and control states},
\newblock \bibinfo{journal}{Infectious Disease Modelling} \bibinfo{volume}{8} (\bibinfo{year}{2023}) \bibinfo{pages}{1015--1031}. \DOIprefix\doi{10.1016/j.idm.2023.08.005}.
\bibitem[{Verma et~al.(2020)Verma, Saini, Gandhi, Dash, and Koya}]{verma2020capacity}
\bibinfo{author}{V.~R. Verma}, \bibinfo{author}{A.~Saini}, \bibinfo{author}{S.~Gandhi}, \bibinfo{author}{U.~Dash}, \bibinfo{author}{S.~F. Koya},
\newblock \bibinfo{title}{{Capacity-need gap in hospital resources for varying mitigation and containment strategies in India in the face of COVID-19 pandemic}},
\newblock \bibinfo{journal}{Infectious Disease Modelling} \bibinfo{volume}{5} (\bibinfo{year}{2020}) \bibinfo{pages}{608--621}. \DOIprefix\doi{10.1016/j.idm.2020.08.011}.
\bibitem[{Zakary et~al.(2017)Zakary, Rachik, and Elmouki}]{zakary2017new}
\bibinfo{author}{O.~Zakary}, \bibinfo{author}{M.~Rachik}, \bibinfo{author}{I.~Elmouki},
\newblock \bibinfo{title}{{A new epidemic modeling approach: Multi-regions discrete-time model with travel-blocking vicinity optimal control strategy}},
\newblock \bibinfo{journal}{Infectious Disease Modelling} \bibinfo{volume}{2} (\bibinfo{year}{2017}) \bibinfo{pages}{304--322}. \DOIprefix\doi{10.1016/j.idm.2017.06.003}.
\bibitem[{Bhadauria et~al.(2023)Bhadauria, Dhungana, Verma, Woodcock, and Rai}]{bhadauria2023studying}
\bibinfo{author}{A.~S. Bhadauria}, \bibinfo{author}{H.~N. Dhungana}, \bibinfo{author}{V.~Verma}, \bibinfo{author}{S.~Woodcock}, \bibinfo{author}{T.~Rai},
\newblock \bibinfo{title}{{Studying the efficacy of isolation as a control strategy and elimination of tuberculosis in India: A mathematical model}},
\newblock \bibinfo{journal}{Infectious Disease Modelling} \bibinfo{volume}{8} (\bibinfo{year}{2023}) \bibinfo{pages}{458--470}. \DOIprefix\doi{10.1016/j.idm.2023.03.005}.
\bibitem[{Vatcheva et~al.(2021)Vatcheva, Sifuentes, Oraby, Maldonado, Huber, and Villalobos}]{vatcheva2021social}
\bibinfo{author}{K.~P. Vatcheva}, \bibinfo{author}{J.~Sifuentes}, \bibinfo{author}{T.~Oraby}, \bibinfo{author}{J.~C. Maldonado}, \bibinfo{author}{T.~Huber}, \bibinfo{author}{M.~C. Villalobos},
\newblock \bibinfo{title}{{Social distancing and testing as optimal strategies against the spread of COVID-19 in the Rio Grande Valley of Texas}},
\newblock \bibinfo{journal}{Infectious Disease Modelling} \bibinfo{volume}{6} (\bibinfo{year}{2021}) \bibinfo{pages}{729--742}. \DOIprefix\doi{10.1016/j.idm.2021.04.004}.
\bibitem[{Brethouwer et~al.(2021)Brethouwer, van~de Rijt, Lindelauf, and Fokkink}]{brethouwer2021stay}
\bibinfo{author}{J.-T. Brethouwer}, \bibinfo{author}{A.~van~de Rijt}, \bibinfo{author}{R.~Lindelauf}, \bibinfo{author}{R.~Fokkink},
\newblock \bibinfo{title}{{“Stay nearby or get checked”: A Covid-19 control strategy}},
\newblock \bibinfo{journal}{Infectious Disease Modelling} \bibinfo{volume}{6} (\bibinfo{year}{2021}) \bibinfo{pages}{36--45}. \DOIprefix\doi{10.1016/j.idm.2020.10.013}.
\bibitem[{Amaku et~al.(2021)Amaku, Covas, Coutinho, Neto, Struchiner, Wilder-Smith, and Massad}]{amaku2021modelling}
\bibinfo{author}{M.~Amaku}, \bibinfo{author}{D.~T. Covas}, \bibinfo{author}{F.~A.~B. Coutinho}, \bibinfo{author}{R.~S.~A. Neto}, \bibinfo{author}{C.~Struchiner}, \bibinfo{author}{A.~Wilder-Smith}, \bibinfo{author}{E.~Massad},
\newblock \bibinfo{title}{{Modelling the test, trace and quarantine strategy to control the COVID-19 epidemic in the state of S{\~a}o Paulo, Brazil}},
\newblock \bibinfo{journal}{Infectious Disease Modelling} \bibinfo{volume}{6} (\bibinfo{year}{2021}) \bibinfo{pages}{46--55}. \DOIprefix\doi{10.1016/j.idm.2020.11.004}.
\bibitem[{Saha et~al.(2022)Saha, Saha, and Podder}]{saha2022effect}
\bibinfo{author}{A.~K. Saha}, \bibinfo{author}{S.~Saha}, \bibinfo{author}{C.~N. Podder},
\newblock \bibinfo{title}{{Effect of awareness, quarantine and vaccination as control strategies on Covid-19 with co-morbidity and re-infection}},
\newblock \bibinfo{journal}{Infectious Disease Modelling} \bibinfo{volume}{7} (\bibinfo{year}{2022}) \bibinfo{pages}{660--689}. \DOIprefix\doi{10.2139/ssrn.4185141}.
\bibitem[{Pat{\'o}n et~al.(2023)Pat{\'o}n, Acu{\~n}a, and Rodr{\'\i}guez}]{paton2023evaluation}
\bibinfo{author}{M.~Pat{\'o}n}, \bibinfo{author}{J.~M. Acu{\~n}a}, \bibinfo{author}{J.~Rodr{\'\i}guez},
\newblock \bibinfo{title}{{Evaluation of vaccine rollout strategies for emerging infectious diseases: A model-based approach including protection attitudes}},
\newblock \bibinfo{journal}{Infectious Disease Modelling} \bibinfo{volume}{8} (\bibinfo{year}{2023}) \bibinfo{pages}{1032--1049}. \DOIprefix\doi{10.1016/j.idm.2023.07.012}.
\bibitem[{Anupong et~al.(2023)Anupong, Chantanasaro, Wilasang, Jitsuk, Sararat, Sornbundit, Pattanasiri, Wannigama, Amarasiri, Chadsuthi, and Modchang}]{anupong2023modeling}
\bibinfo{author}{S.~Anupong}, \bibinfo{author}{T.~Chantanasaro}, \bibinfo{author}{C.~Wilasang}, \bibinfo{author}{N.~C. Jitsuk}, \bibinfo{author}{C.~Sararat}, \bibinfo{author}{K.~Sornbundit}, \bibinfo{author}{B.~Pattanasiri}, \bibinfo{author}{D.~L. Wannigama}, \bibinfo{author}{M.~Amarasiri}, \bibinfo{author}{S.~Chadsuthi}, \bibinfo{author}{C.~Modchang},
\newblock \bibinfo{title}{{Modeling vaccination strategies with limited early {COVID}-19 vaccine access in low-and middle-income countries: A case study of Thailand}},
\newblock \bibinfo{journal}{Infectious Disease Modelling} \bibinfo{volume}{8} (\bibinfo{year}{2023}) \bibinfo{pages}{1177--1189}. \DOIprefix\doi{10.1016/j.idm.2023.11.003}.
\bibitem[{Owusu-Dampare and Bouchnita(2023)}]{owusu2023equitable}
\bibinfo{author}{F.~Owusu-Dampare}, \bibinfo{author}{A.~Bouchnita},
\newblock \bibinfo{title}{{Equitable bivalent booster allocation strategies against emerging SARS-CoV-2 variants in US cities with large Hispanic communities: The case of El Paso County, Texas}},
\newblock \bibinfo{journal}{Infectious Disease Modelling} \bibinfo{volume}{8} (\bibinfo{year}{2023}) \bibinfo{pages}{912--919}. \DOIprefix\doi{10.1016/j.idm.2023.07.009}.
\bibitem[{Gan et~al.(2024)Gan, Janhavi, Tong, Lim, and Dickens}]{gan2024need}
\bibinfo{author}{G.~Gan}, \bibinfo{author}{A.~Janhavi}, \bibinfo{author}{G.~Tong}, \bibinfo{author}{J.~T. Lim}, \bibinfo{author}{B.~L. Dickens},
\newblock \bibinfo{title}{{The need for pre-emptive control strategies for mpox in Asia and Oceania}},
\newblock \bibinfo{journal}{Infectious Disease Modelling} \bibinfo{volume}{9} (\bibinfo{year}{2024}) \bibinfo{pages}{214--223}. \DOIprefix\doi{10.1016/j.idm.2023.12.005}.
\bibitem[{Thongtha and Modnak(2022)}]{thongtha2022optimal}
\bibinfo{author}{A.~Thongtha}, \bibinfo{author}{C.~Modnak},
\newblock \bibinfo{title}{{Optimal COVID-19 epidemic strategy with vaccination control and infection prevention measures in Thailand}},
\newblock \bibinfo{journal}{Infectious Disease Modelling} \bibinfo{volume}{7} (\bibinfo{year}{2022}) \bibinfo{pages}{835--855}. \DOIprefix\doi{10.1016/j.idm.2022.11.002}.
\bibitem[{Abell et~al.(2023)Abell, McCaw, and Baker}]{abell2023understanding}
\bibinfo{author}{I.~R. Abell}, \bibinfo{author}{J.~M. McCaw}, \bibinfo{author}{C.~M. Baker},
\newblock \bibinfo{title}{Understanding the impact of disease and vaccine mechanisms on the importance of optimal vaccine allocation},
\newblock \bibinfo{journal}{Infectious Disease Modelling} \bibinfo{volume}{8} (\bibinfo{year}{2023}) \bibinfo{pages}{539--550}. \DOIprefix\doi{10.1016/j.idm.2023.05.003}.
\bibitem[{Zaman et~al.(2009)Zaman, Kang, and Jung}]{ZAMAN200943}
\bibinfo{author}{G.~Zaman}, \bibinfo{author}{Y.~H. Kang}, \bibinfo{author}{I.~H. Jung},
\newblock \bibinfo{title}{{Optimal treatment of an SIR epidemic model with time delay}},
\newblock \bibinfo{journal}{Biosystems} \bibinfo{volume}{98} (\bibinfo{year}{2009}) \bibinfo{pages}{43--50}. \DOIprefix\doi{doi.org/10.1016/j.biosystems.2009.05.006}.
\bibitem[{B{\'e}raud et~al.(2022)B{\'e}raud, Timsit, and Leleu}]{beraud2022remdesivir}
\bibinfo{author}{G.~B{\'e}raud}, \bibinfo{author}{J.-F. Timsit}, \bibinfo{author}{H.~Leleu},
\newblock \bibinfo{title}{{Remdesivir and dexamethasone as tools to relieve hospital care systems stressed by COVID-19: a modelling study on bed resources and budget impact}},
\newblock \bibinfo{journal}{PLOS ONE} \bibinfo{volume}{17} (\bibinfo{year}{2022}) \bibinfo{pages}{e0262462}. \DOIprefix\doi{10.1371/journal.pone.0262462}.
\bibitem[{Apenteng et~al.(2020)Apenteng, Osei, Oduro, Kwabla, and Ismail}]{apenteng2020impact}
\bibinfo{author}{O.~O. Apenteng}, \bibinfo{author}{P.~P. Osei}, \bibinfo{author}{B.~Oduro}, \bibinfo{author}{M.~P. Kwabla}, \bibinfo{author}{N.~A. Ismail},
\newblock \bibinfo{title}{{The impact of implementing HIV prevention policies therapy and control strategy among HIV and AIDS incidence cases in Malaysia}},
\newblock \bibinfo{journal}{Infectious Disease Modelling} \bibinfo{volume}{5} (\bibinfo{year}{2020}) \bibinfo{pages}{755--765}. \DOIprefix\doi{10.1016/j.idm.2020.09.009}.
\bibitem[{Lamba et~al.(2024)Lamba, Das, and Srivastava}]{lamba2024impact}
\bibinfo{author}{S.~Lamba}, \bibinfo{author}{T.~Das}, \bibinfo{author}{P.~K. Srivastava},
\newblock \bibinfo{title}{{Impact of infectious density-induced additional screening and treatment saturation on COVID-19: Modeling and cost-effective optimal control}},
\newblock \bibinfo{journal}{Infectious Disease Modelling} \bibinfo{volume}{9} (\bibinfo{year}{2024}) \bibinfo{pages}{569--600}. \DOIprefix\doi{10.1016/j.idm.2024.03.002}.
\bibitem[{Van~Rens and Oswald(2020)}]{van2020age}
\bibinfo{author}{T.~Van~Rens}, \bibinfo{author}{A.~J. Oswald},
\newblock \bibinfo{title}{Age-based policy in the context of the {Covid}-19 pandemic: how common are multigenerational households?},
\newblock \bibinfo{journal}{CAGE Online Working Paper Series}  (\bibinfo{year}{2020}).
\bibitem[{Spaccatini et~al.(2022)Spaccatini, Giovannelli, and Pacilli}]{spaccatini2022you}
\bibinfo{author}{F.~Spaccatini}, \bibinfo{author}{I.~Giovannelli}, \bibinfo{author}{M.~G. Pacilli},
\newblock \bibinfo{title}{“{Y}ou are stealing our present”: {Y}ounger people's ageism towards older people predicts attitude towards age-based {COVID}-19 restriction measures},
\newblock \bibinfo{journal}{Journal of Social Issues} \bibinfo{volume}{78} (\bibinfo{year}{2022}) \bibinfo{pages}{769--789}. \DOIprefix\doi{doi.org/10.1111/josi.12537}.
\bibitem[{Motorniak et~al.(2023)Motorniak, Savulescu, and Giubilini}]{motorniak2023reelin}
\bibinfo{author}{D.~Motorniak}, \bibinfo{author}{J.~Savulescu}, \bibinfo{author}{A.~Giubilini},
\newblock \bibinfo{title}{{R}eelin’in the years: Age and selective restriction of liberty in the {COVID}-19 pandemic},
\newblock \bibinfo{journal}{Journal of Bioethical Inquiry}  (\bibinfo{year}{2023}) \bibinfo{pages}{1--9}. \DOIprefix\doi{10.1007/s11673-023-10318-8}.
\bibitem[{Acemoglu et~al.(2021)Acemoglu, Chernozhukov, Werning, and Whinston}]{acemoglu2021optimal}
\bibinfo{author}{D.~Acemoglu}, \bibinfo{author}{V.~Chernozhukov}, \bibinfo{author}{I.~Werning}, \bibinfo{author}{M.~D. Whinston},
\newblock \bibinfo{title}{Optimal targeted lockdowns in a multigroup sir model},
\newblock \bibinfo{journal}{American Economic Review: Insights} \bibinfo{volume}{3} (\bibinfo{year}{2021}) \bibinfo{pages}{487--502}. \DOIprefix\doi{10.1257/aeri.20200590}.
\bibitem[{Kirwin et~al.(2021)Kirwin, Rafferty, Harback, Round, and McCabe}]{kirwin2021net}
\bibinfo{author}{E.~Kirwin}, \bibinfo{author}{E.~Rafferty}, \bibinfo{author}{K.~Harback}, \bibinfo{author}{J.~Round}, \bibinfo{author}{C.~McCabe},
\newblock \bibinfo{title}{A net benefit approach for the optimal allocation of a covid-19 vaccine},
\newblock \bibinfo{journal}{Pharmacoeconomics} \bibinfo{volume}{39} (\bibinfo{year}{2021}) \bibinfo{pages}{1059--1073}. \DOIprefix\doi{10.1007/s40273-021-01037-2}.
\bibitem[{Diekmann and Heesterbeek(2000)}]{diekmann2000mathematical}
\bibinfo{author}{O.~Diekmann}, \bibinfo{author}{J.~A.~P. Heesterbeek}, \bibinfo{title}{{Mathematical Epidemiology of Infectious Diseases: Model Building, Analysis and Interpretation}}, volume~\bibinfo{volume}{5}, \bibinfo{publisher}{John Wiley \& Sons}, \bibinfo{year}{2000}.
\bibitem[{Del~Valle et~al.(2007)Del~Valle, Hyman, Hethcote, and Eubank}]{DelValle_Hyman_Hethcote_Eubank_2007}
\bibinfo{author}{S.~Del~Valle}, \bibinfo{author}{J.~Hyman}, \bibinfo{author}{H.~Hethcote}, \bibinfo{author}{S.~Eubank},
\newblock \bibinfo{title}{Mixing patterns between age groups in social networks},
\newblock \bibinfo{journal}{Social Networks} \bibinfo{volume}{29} (\bibinfo{year}{2007}) \bibinfo{pages}{539–554}. \DOIprefix\doi{10.1016/j.socnet.2007.04.005}.
\bibitem[{Grant et~al.(2022)Grant, Charmet, Schaeffer, Galmiche, Madec, Von~Platen, Ch{\'e}ny, Omar, David, Rogoff et~al.}]{grant2022impact}
\bibinfo{author}{R.~Grant}, \bibinfo{author}{T.~Charmet}, \bibinfo{author}{L.~Schaeffer}, \bibinfo{author}{S.~Galmiche}, \bibinfo{author}{Y.~Madec}, \bibinfo{author}{C.~Von~Platen}, \bibinfo{author}{O.~Ch{\'e}ny}, \bibinfo{author}{F.~Omar}, \bibinfo{author}{C.~David}, \bibinfo{author}{A.~Rogoff}, et~al.,
\newblock \bibinfo{title}{{Impact of SARS-CoV-2 Delta variant on incubation, transmission settings and vaccine effectiveness: Results from a nationwide case-control study in France}},
\newblock \bibinfo{journal}{The Lancet Regional Health - Europe} \bibinfo{volume}{13} (\bibinfo{year}{2022}) \bibinfo{pages}{100278}. \DOIprefix\doi{10.1016/j.lanepe.2021.100278}.
\bibitem[{Chau et~al.(2022)Chau, Nguyet, Truong, Dung, Nhan, Man, Ngoc, Thao, Tu, Mai et~al.}]{chau2022immunogenicity}
\bibinfo{author}{N.~V.~V. Chau}, \bibinfo{author}{L.~A. Nguyet}, \bibinfo{author}{N.~T. Truong}, \bibinfo{author}{N.~T. Dung}, \bibinfo{author}{M.~T. Nhan}, \bibinfo{author}{D.~N.~H. Man}, \bibinfo{author}{N.~M. Ngoc}, \bibinfo{author}{H.~P. Thao}, \bibinfo{author}{T.~N.~H. Tu}, \bibinfo{author}{H.~K. Mai}, et~al.,
\newblock \bibinfo{title}{{Immunogenicity of Oxford-AstraZeneca COVID-19 vaccine in Vietnamese health-care workers}},
\newblock \bibinfo{journal}{The American Journal of Tropical Medicine and Hygiene} \bibinfo{volume}{106} (\bibinfo{year}{2022}) \bibinfo{pages}{556}. \DOIprefix\doi{10.4269/ajtmh.21-0849}.
\bibitem[{Kang et~al.(2022)Kang, Xin, Yuan, Ali, Liang, Zhang, Hu, Lau, Zhang, Zhang et~al.}]{kang2022transmission}
\bibinfo{author}{M.~Kang}, \bibinfo{author}{H.~Xin}, \bibinfo{author}{J.~Yuan}, \bibinfo{author}{S.~T. Ali}, \bibinfo{author}{Z.~Liang}, \bibinfo{author}{J.~Zhang}, \bibinfo{author}{T.~Hu}, \bibinfo{author}{E.~H. Lau}, \bibinfo{author}{Y.~Zhang}, \bibinfo{author}{M.~Zhang}, et~al.,
\newblock \bibinfo{title}{{Transmission dynamics and epidemiological characteristics of SARS-CoV-2 Delta variant infections in Guangdong, China, May to June 2021}},
\newblock \bibinfo{journal}{Eurosurveillance} \bibinfo{volume}{27} (\bibinfo{year}{2022}) \bibinfo{pages}{2100815}. \DOIprefix\doi{10.2807/1560-7917.ES.2022.27.10.2100815}.
\bibitem[{Wu et~al.(2022)Wu, Kang, Guo, Liu, Liu, and Liang}]{wu2022incubation}
\bibinfo{author}{Y.~Wu}, \bibinfo{author}{L.~Kang}, \bibinfo{author}{Z.~Guo}, \bibinfo{author}{J.~Liu}, \bibinfo{author}{M.~Liu}, \bibinfo{author}{W.~Liang},
\newblock \bibinfo{title}{{Incubation period of COVID-19 caused by unique SARS-CoV-2 strains: a systematic review and meta-analysis}},
\newblock \bibinfo{journal}{JAMA Network Open} \bibinfo{volume}{5} (\bibinfo{year}{2022}) \bibinfo{pages}{e2228008--e2228008}. \DOIprefix\doi{10.1001/jamanetworkopen.2022.28008}.
\bibitem[{Sah et~al.(2021)Sah, Fitzpatrick, Zimmer, Abdollahi, Juden-Kelly, Moghadas, Singer, and Galvani}]{sah2021asymptomatic}
\bibinfo{author}{P.~Sah}, \bibinfo{author}{M.~C. Fitzpatrick}, \bibinfo{author}{C.~F. Zimmer}, \bibinfo{author}{E.~Abdollahi}, \bibinfo{author}{L.~Juden-Kelly}, \bibinfo{author}{S.~M. Moghadas}, \bibinfo{author}{B.~H. Singer}, \bibinfo{author}{A.~P. Galvani},
\newblock \bibinfo{title}{{Asymptomatic SARS-CoV-2 infection: A systematic review and meta-analysis}},
\newblock \bibinfo{journal}{Proceedings of the National Academy of Sciences} \bibinfo{volume}{118} (\bibinfo{year}{2021}) \bibinfo{pages}{e2109229118}. \DOIprefix\doi{10.1073/pnas.2109229118}.
\bibitem[{He et~al.(2021)He, Guo, Mao, and Zhang}]{he2021proportion}
\bibinfo{author}{J.~He}, \bibinfo{author}{Y.~Guo}, \bibinfo{author}{R.~Mao}, \bibinfo{author}{J.~Zhang},
\newblock \bibinfo{title}{Proportion of asymptomatic coronavirus disease 2019: A systematic review and meta-analysis},
\newblock \bibinfo{journal}{Journal of Medical Virology} \bibinfo{volume}{93} (\bibinfo{year}{2021}) \bibinfo{pages}{820--830}. \DOIprefix\doi{10.1002/jmv.26326}.
\bibitem[{Buitrago-Garcia et~al.(2022)Buitrago-Garcia, Ipekci, Heron, Imeri, Araujo-Chaveron, Arevalo-Rodriguez, Ciapponi, Cevik, Hauser, Alam et~al.}]{buitrago2022occurrence}
\bibinfo{author}{D.~Buitrago-Garcia}, \bibinfo{author}{A.~M. Ipekci}, \bibinfo{author}{L.~Heron}, \bibinfo{author}{H.~Imeri}, \bibinfo{author}{L.~Araujo-Chaveron}, \bibinfo{author}{I.~Arevalo-Rodriguez}, \bibinfo{author}{A.~Ciapponi}, \bibinfo{author}{M.~Cevik}, \bibinfo{author}{A.~Hauser}, \bibinfo{author}{M.~I. Alam}, et~al.,
\newblock \bibinfo{title}{{Occurrence and transmission potential of asymptomatic and presymptomatic SARS-CoV-2 infections: Update of a living systematic review and meta-analysis}},
\newblock \bibinfo{journal}{PLOS Medicine} \bibinfo{volume}{19} (\bibinfo{year}{2022}) \bibinfo{pages}{e1003987}. \DOIprefix\doi{10.1371/journal.pmed.1003987}.
\bibitem[{Byrne et~al.(2020)Byrne, McEvoy, Collins, Hunt, Casey, Barber, Butler, Griffin, Lane, McAloon et~al.}]{byrne2020inferred}
\bibinfo{author}{A.~W. Byrne}, \bibinfo{author}{D.~McEvoy}, \bibinfo{author}{A.~B. Collins}, \bibinfo{author}{K.~Hunt}, \bibinfo{author}{M.~Casey}, \bibinfo{author}{A.~Barber}, \bibinfo{author}{F.~Butler}, \bibinfo{author}{J.~Griffin}, \bibinfo{author}{E.~A. Lane}, \bibinfo{author}{C.~McAloon}, et~al.,
\newblock \bibinfo{title}{Inferred duration of infectious period of {SARS}-{C}o{V}-2: rapid scoping review and analysis of available evidence for asymptomatic and symptomatic {COVID}-19 cases},
\newblock \bibinfo{journal}{BMJ Open} \bibinfo{volume}{10} (\bibinfo{year}{2020}) \bibinfo{pages}{e039856}. \DOIprefix\doi{10.1136/bmjopen-2020-039856}.
\bibitem[{{Centers for Disease Control and Prevention}(2020)}]{CDCP2020}
\bibinfo{author}{{Centers for Disease Control and Prevention}}, \bibinfo{year}{2020}. \URLprefix \url{https://www.cdc.gov/coronavirus/2019-ncov/hcp/testing-overview.html}.
\bibitem[{Pollock and Lancaster(2020)}]{pollock2020asymptomatic}
\bibinfo{author}{A.~M. Pollock}, \bibinfo{author}{J.~Lancaster},
\newblock \bibinfo{title}{Asymptomatic transmission of {C}ovid-19},
\newblock \bibinfo{journal}{BMJ} \bibinfo{volume}{371} (\bibinfo{year}{2020}). \DOIprefix\doi{10.1136/bmj.m4851}.
\bibitem[{{SAGE 56th meeting on COVID-19}(2020)}]{sage56}
\bibinfo{author}{{SAGE 56th meeting on COVID-19}}, \bibinfo{year}{2020}. \bibinfo{note}{\url{https://assets.publishing.service.gov.uk/government/uploads/system/uploads/attachment_data/file/928699/S0740_Fifty-sixth_SAGE_meeting_on_Covid-19.pdf}}.
\bibitem[{Soni et~al.(2023)Soni, Herbert, Lin, Yan, Pretz, Stamegna, Wang, Orwig, Wright, Tarrant et~al.}]{soni2023performance}
\bibinfo{author}{A.~Soni}, \bibinfo{author}{C.~Herbert}, \bibinfo{author}{H.~Lin}, \bibinfo{author}{Y.~Yan}, \bibinfo{author}{C.~Pretz}, \bibinfo{author}{P.~Stamegna}, \bibinfo{author}{B.~Wang}, \bibinfo{author}{T.~Orwig}, \bibinfo{author}{C.~Wright}, \bibinfo{author}{S.~Tarrant}, et~al.,
\newblock \bibinfo{title}{{Performance of Rapid Antigen Tests to Detect Symptomatic and Asymptomatic SARS-CoV-2 Infection: A Prospective Cohort Study}},
\newblock \bibinfo{journal}{Annals of Internal Medicine} \bibinfo{volume}{176} (\bibinfo{year}{2023}) \bibinfo{pages}{975--982}. \DOIprefix\doi{10.7326/M23-0385}.
\bibitem[{Billah et~al.(2020)Billah, Miah, and Khan}]{billah2020reproductive}
\bibinfo{author}{M.~A. Billah}, \bibinfo{author}{M.~M. Miah}, \bibinfo{author}{M.~N. Khan},
\newblock \bibinfo{title}{Reproductive number of coronavirus: A systematic review and meta-analysis based on global level evidence},
\newblock \bibinfo{journal}{PLOS ONE} \bibinfo{volume}{15} (\bibinfo{year}{2020}) \bibinfo{pages}{e0242128}. \DOIprefix\doi{10.1371/journal.pone.0242128}.
\bibitem[{Ahammed et~al.(2021)Ahammed, Anjum, Rahman, Haider, Kock, and Uddin}]{ahammed2021estimation}
\bibinfo{author}{T.~Ahammed}, \bibinfo{author}{A.~Anjum}, \bibinfo{author}{M.~M. Rahman}, \bibinfo{author}{N.~Haider}, \bibinfo{author}{R.~Kock}, \bibinfo{author}{M.~J. Uddin},
\newblock \bibinfo{title}{Estimation of novel coronavirus ({COVID}-19) reproduction number and case fatality rate: A systematic review and meta-analysis},
\newblock \bibinfo{journal}{Health Science Reports} \bibinfo{volume}{4} (\bibinfo{year}{2021}) \bibinfo{pages}{e274}. \DOIprefix\doi{10.1002/hsr2.274}.
\bibitem[{Vinceti et~al.(2022)Vinceti, Balboni, Rothman, Teggi, Bellino, Pezzotti, Ferrari, Orsini, and Filippini}]{vinceti2022substantial}
\bibinfo{author}{M.~Vinceti}, \bibinfo{author}{E.~Balboni}, \bibinfo{author}{K.~J. Rothman}, \bibinfo{author}{S.~Teggi}, \bibinfo{author}{S.~Bellino}, \bibinfo{author}{P.~Pezzotti}, \bibinfo{author}{F.~Ferrari}, \bibinfo{author}{N.~Orsini}, \bibinfo{author}{T.~Filippini},
\newblock \bibinfo{title}{{Substantial impact of mobility restrictions on reducing {COVID}-19 incidence in {I}taly in 2020}},
\newblock \bibinfo{journal}{Journal of Travel Medicine} \bibinfo{volume}{29} (\bibinfo{year}{2022}). \DOIprefix\doi{10.1093/jtm/taac081}.
\bibitem[{Schlosser et~al.(2020)Schlosser, Maier, Jack, Hinrichs, Zachariae, and Brockmann}]{schlosser2020covid}
\bibinfo{author}{F.~Schlosser}, \bibinfo{author}{B.~F. Maier}, \bibinfo{author}{O.~Jack}, \bibinfo{author}{D.~Hinrichs}, \bibinfo{author}{A.~Zachariae}, \bibinfo{author}{D.~Brockmann},
\newblock \bibinfo{title}{{COVID}-19 lockdown induces disease-mitigating structural changes in mobility networks},
\newblock \bibinfo{journal}{Proceedings of the National Academy of Sciences} \bibinfo{volume}{117} (\bibinfo{year}{2020}) \bibinfo{pages}{32883--32890}. \DOIprefix\doi{10.1073/pnas.2012326117}.
\bibitem[{Zhou et~al.(2020)Zhou, Xu, Hu, Yue, Li, and Xia}]{zhou2020effects}
\bibinfo{author}{Y.~Zhou}, \bibinfo{author}{R.~Xu}, \bibinfo{author}{D.~Hu}, \bibinfo{author}{Y.~Yue}, \bibinfo{author}{Q.~Li}, \bibinfo{author}{J.~Xia},
\newblock \bibinfo{title}{Effects of human mobility restrictions on the spread of {COVID}-19 in {S}henzhen, {C}hina: a modelling study using mobile phone data},
\newblock \bibinfo{journal}{The Lancet Digital Health} \bibinfo{volume}{2} (\bibinfo{year}{2020}) \bibinfo{pages}{e417--e424}. \DOIprefix\doi{10.1016/S2589-7500(20)30165-5}.
\bibitem[{{Wolfram Research Inc.}(2022)}]{Mathematica}
\bibinfo{author}{{Wolfram Research Inc.}}, \bibinfo{title}{Mathematica, {V}ersion 13.1}, \bibinfo{year}{2022}. \bibinfo{note}{\url{https://www.wolfram.com/mathematica}}.
\bibitem[{Gao et~al.(2021)Gao, Xu, Sun, Wang, Guo, Qiu, and Ma}]{gao2021systematic}
\bibinfo{author}{Z.~Gao}, \bibinfo{author}{Y.~Xu}, \bibinfo{author}{C.~Sun}, \bibinfo{author}{X.~Wang}, \bibinfo{author}{Y.~Guo}, \bibinfo{author}{S.~Qiu}, \bibinfo{author}{K.~Ma},
\newblock \bibinfo{title}{{A systematic review of asymptomatic infections with COVID-19}},
\newblock \bibinfo{journal}{Journal of Microbiology, Immunology and Infection} \bibinfo{volume}{54} (\bibinfo{year}{2021}) \bibinfo{pages}{12--16}. \DOIprefix\doi{10.1016/j.jmii.2020.05.001}.
\bibitem[{Rohatgi(2022)}]{Rohatgi2022}
\bibinfo{author}{A.~Rohatgi}, \bibinfo{title}{Webplotdigitizer: Version 4.6}, \bibinfo{howpublished}{\url{https://automeris.io/WebPlotDigitizer}}, \bibinfo{year}{2022}.
\bibitem[{Tan et~al.(2020)Tan, Wong, Leo, and Toh}]{tan2020does}
\bibinfo{author}{W.~Tan}, \bibinfo{author}{L.~Wong}, \bibinfo{author}{Y.~Leo}, \bibinfo{author}{M.~Toh},
\newblock \bibinfo{title}{{Does incubation period of COVID-19 vary with age? A study of epidemiologically linked cases in Singapore}},
\newblock \bibinfo{journal}{Epidemiology \& Infection} \bibinfo{volume}{148} (\bibinfo{year}{2020}) \bibinfo{pages}{e197}. \DOIprefix\doi{10.1017/S0950268820001995}.

\end{thebibliography}

\end{document}